\documentclass{compositio}

\usepackage{rotating,amsthm,amsfonts,amsmath,amssymb,mathrsfs,verbatim,cite}

\numberwithin{equation}{section}

\newtheorem{lemma}{Lemma}[section]
\newtheorem{proposition}{Proposition}[section]
\newtheorem{theorem}{Theorem}[section]
\newtheorem{corollary}{Corollary}[section]
\newtheorem{definition}{Definition}[section]

\renewcommand{\Re}{\textup{Re}}
\newcommand{\dup}{\textup{d}}

\newcommand{\Int}{\int\limits}
\newcommand{\abs}[1]{\lvert#1\rvert}
\newcommand{\F}{\mathbb F}
\newcommand{\Q}{\mathbb Q}

\newcommand{\Z}{\mathbb Z}
\newcommand{\Complex}{\mathbb C}

\newcommand{\la}{\lambda}

\newcommand{\gsl}{\mathfrak{sl}}
\newcommand{\Symm}{\mathfrak{S}}
\newcommand{\qbin}[2]{\genfrac{[}{]}{0pt}{}{#1}{#2}}
\newcommand{\hyp}[3]{\biggl[\genfrac{}{}{0pt}{}{#1}{#2};#3\biggr]}
\newcommand{\LEQ}{\!\!\leq\!\!}
\newcommand{\GEQ}{\!\!\geq\!\!}

\begin{document}

\title[$\gsl_3$ Basic hypergeometric series]
{Bisymmetric functions, Macdonald polynomials and \\
$\gsl_3$ basic hypergeometric series}

\author{S. Ole Warnaar}\thanks{Work supported by the Australian 
Research Council}

\email{warnaar@ms.unimelb.edu.au}

\address{Department of Mathematics and Statistics,
The University of Melbourne, VIC 3010, Australia}



\begin{abstract}
A new type of $\gsl_3$ basic hypergeometric series 
based on Macdonald polynomials is introduced.
Besides a pair of Macdonald polynomials attached to two different
sets of variables, a key-ingredient in the
$\gsl_3$ basic hypergeometric series is a bisymmetric function related to
Macdonald's commuting family of $q$-difference operators, 
to the $\gsl_3$ Selberg integrals of Tarasov and Varchenko,
and to alternating sign matrices.
Our main result for $\gsl_3$ series is a multivariable generalization of the 
celebrated $q$-binomial theorem. In the limit this $q$-binomial sum
yields a new $\gsl_3$ Selberg integral for Jack polynomials.
\end{abstract}

\maketitle

\section{Introduction}
The $q$-binomial theorem, which was independently discovered by 
Cauchy, Heine and Gauss (with special cases due to Euler and Rothe) 
is one of the most important results in the theory of $q$-series,
see e.g., \cite{AAR99,GR04} and references therein.
Using the standard notation $(a;q)_n=(1-a)(1-aq)\cdots(1-aq^{n-1})$ 
for the $q$-shifted factorial, the theorem may be stated as
\begin{equation}\label{qbthm}
{_1\phi_0}\hyp{a}{\text{--}}{q,z}:=
\sum_{k=0}^{\infty} \frac{(a;q)_k}{(q;q)_k}\, z^k
=\frac{(az;q)_{\infty}}{(z;q)_{\infty}}
\end{equation}
for $\abs{q}<1$ and $\abs{z}<1$.
A well-known alternative representation of the $q$-binomial theorem
is as the $q$-beta integral (for the definition of $q$-integrals
see \cite{GR04})
\[
\int_0^1 t^{\alpha-1}(tq;q)_{\beta-1}\,
\dup_q t=
\frac{\Gamma_q(\alpha)\Gamma_q(\beta)}{\Gamma_q(\alpha+\beta)},
\]
where $0<q<1$, $\Gamma_q$ is the $q$-gamma function \cite{GR04}, 
\[
(a;q)_z=\frac{(a;q)_{\infty}}{(aq^z;q)_{\infty}}\quad\text{for $z\in\Complex$},
\]
and $\alpha,\beta\in\Complex$ such that
$\Re(\alpha)>0$, $-\beta\not\in\{0,1,2,\dots\}$.
Assuming $\Re(\beta)>0$ and taking the limit $q\to 1^{-}$ it follows
that the $q$-binomial theorem implies Euler's beta integral \cite{AAR99}
\[
\int_0^1 t^{\alpha-1}(1-t)^{\beta-1}\, \dup t=
\frac{\Gamma(\alpha)\Gamma(\beta)}{\Gamma(\alpha+\beta)}.
\]

Building on the pioneering work of Milne and Gustafson on
multivariable basic hypergeometric series,
many generalizations of the $q$-binomial theorem
have been found in recent times. Most of these
are labelled by one of the classical root systems,
see e.g., \cite{GK96,Milne85,Milne97,ML95,BS98}.
A particularly interesting generalisation of the 
$q$-binomial series is obtained when $z^k$ in \eqref{qbthm}
is replaced by an appropriate symmetric function 
such as the Schur function or Macdonald polynomial,
see \cite{BF99,Kaneko96,Macdonald,Milne92}.
The latter case was independently 
considered by Kaneko and Macdonald, who proved that
\cite{Kaneko96,Macdonald}
\begin{equation}\label{qbthm2}
{_1\Phi_0}\hyp{a}{\text{--}}{q,t;x}
:=\sum_{\lambda} t^{n(\la)}
\frac{(a;q,t)_{\lambda}}{c'_{\lambda}(q,t)}\:
P_{\lambda}(x;q,t)
=\prod_{i\geq 1} \frac{(ax_i;q)_{\infty}}{(x_i;q)_{\infty}}.
\end{equation}
Here $P_{\la}(x;q,t)$ is the Macdonald
polynomial labelled by the partition $\la$, 
$n(\la)=\sum_{i\geq 1} (i-1)\la_i$, and $c'_{\la}(q,t)$ and
$(a;q,t)_{\la}$ (defined in Section~\ref{sec21}) 
are generalisations of the $q$-shifted factorials
$(q;q)_k$ and $(a;q)_k$, respectively.
If $x$ contains a single variable then the partition
$\la$ is restricted to only one part, and
\eqref{qbthm2} reduces to the ordinary $q$-binomial theorem
\eqref{qbthm}.

Analogous to the single-variable case, \eqref{qbthm2} may
be transformed into a multiple $q$-integral.
In the $q\to 1^{-}$ limit this implies the famous
Selberg integral \cite{Selberg44}
\begin{equation}\label{SelbergInt}
\Int_{[0,1]^n}
\prod_{i=1}^n x_i^{\alpha-1}(1-x_i)^{\beta-1}
\prod_{1\leq i<j\leq n} \abs{x_i-x_j}^{2\gamma}
\: \dup x 
=\prod_{i=1}^n
\frac{\Gamma(\alpha+(i-1)\gamma)\Gamma(\beta+(i-1)\gamma)\Gamma(i\gamma+1)}
{\Gamma(\alpha+\beta+(n+i-2)\gamma)\Gamma(\gamma+1)}
\end{equation}
for $\Re(\alpha)>0,~\Re(\beta)>0,~\Re(\gamma)>
-\min\{1/n,\Re(\alpha)/(n-1),\Re(\beta)/(n-1)\}$.

In this paper we take the natural next step in the
development of basic hypergeometric series and
prove an $\gsl_3$ version of the Kaneko--Macdonald
$q$-binomial theorem:
\begin{equation}\label{qbthmsl3}
{_1\Phi_0}\hyp{a}{\textbf{--}}{q,t;x,y}
=\prod_{i=1}^m \frac{(azt^{m-1}x_i;q)_{\infty}}
{(zt^{m-1}x_i;q)_{\infty}}
\prod_{i=1}^{n-m} \frac{(azt^{n-i};q)_{\infty}}{(zt^{n-i};q)_{\infty}}
\end{equation}
for $y=z(1,t,\dots,t^{n-1})$ and $0\leq m\leq n$.
The series on the left (defined in Section~\ref{SecA2}) depends on two
Macdonald polynomials, $P_{\la}(x_1,\dots,x_m;q,t)$ and 
$P_{\mu}(y_1,\dots,y_n;q,t)$,
and --- as a new ingredient --- involves a bisymmetric function
related to Macdonald's commuting family of $q$-difference operators
\cite{Macdonald95}.

As in the previous two cases one may transform the $\gsl_3$
basic hypergeometric series into a multiple $q$-integral.
The $q\to 1^{-}$ limit then yields the $\gsl_3$ Selberg integral 
of Tarasov and Varchenko \cite{TV03}
\begin{align}\label{integral}
&\Int_{C^{m,n}_{\gamma}[0,1]}
h(x,y)\,
\prod_{i=1}^mx_i^{\beta_1-1} 
\prod_{i=1}^n (1-y_i)^{\alpha-1}y_i^{\beta_2-1} \\ 
&\qquad\quad\times
\prod_{1\leq i<j\leq m} \abs{x_i-x_j}^{2\gamma}
\prod_{1\leq i<j\leq n} \abs{y_i-y_j}^{2\gamma}
\prod_{i=1}^m\prod_{j=1}^n \abs{x_i-y_j}^{-\gamma}
\; \dup x\,\dup y \notag \\[2mm]
&\qquad\qquad =\prod_{i=1}^m \frac{\Gamma(\beta_1+(i-1)\gamma)
\Gamma(\beta_1+\beta_2+(i-2)\gamma)
\Gamma((i-n-1)\gamma)\Gamma(i\gamma)}
{\Gamma(\beta_1+(i+m-n-2)\gamma)
\Gamma(\alpha+\beta_1+\beta_2+(i+n-3)\gamma)
\Gamma(\gamma)} \notag \\
&\qquad\qquad\quad \times
\prod_{i=1}^n \frac{\Gamma(\alpha+(i-1)\gamma)\Gamma(i\gamma)}
{\Gamma(\gamma)}
\prod_{i=1}^{n-m}\frac{\Gamma(\beta_2+(i-1)\gamma)}
{\Gamma(\alpha+\beta_2+(i+n-2)\gamma)},
\notag 
\end{align}
where $C^{m,n}_{\gamma}[0,1]$ is an integration domain described in
Section~\ref{SecA2}, $h(x,y)$ is the bisymmetric function
\[
h(x,y)=\frac{(n-m)!}{n!}\sum_{\substack{l_1,\dots,l_m=1 \\ l_i\neq l_j}}^n
\prod_{i=1}^m \frac{y_{l_i}}{y_{l_i}-x_i}
\]
and (for generic $n$ and $m$)
\begin{gather*}
\Re(\alpha)>0,~\Re(\beta_1)>0,~\Re(\beta_2)>0 \\[2mm]
-\min\Bigl\{
\frac{1}{n},\frac{\Re(\alpha)}{n-1},\frac{\Re(\beta_1)}{m-1},
\frac{\Re(\beta_2)}{n-m-1},\frac{\Re(\beta_1+\beta_2)}{m-2}\Bigr\}
<\Re(\gamma)<0.
\end{gather*}

\subsection{Outline}
In the next section we provide a brief introduction to
Macdonald polynomials and the $\gsl_2$ Kaneko--Macdonald
multivariable basic hypergeometric series.
Then, in Section~\ref{secbiF}, we define the bisymmetric function
$F(x,y;t)$, which plays a key-part in the $\gsl_3$ basic
hypergeometric series studied in this paper.
We prove several elementary results for $F$, and establish a connection
with the bisymmetric function of Tarasov and Varchenko, and with
alternating sign matrices.
In Section~\ref{secqtLR} we obtain an identity involving
the $q,t$-Littlewood--Richardson coefficients and a
specialization of the function $F$.
This identity is at the heart of our proof of the $\gsl_3$ $q$-binomial
theorem \eqref{qbthmsl3}.
Finally, in Section~\ref{SecA2} we define the $\gsl_3$ basic 
hypergeometric series and prove several $q$-binomial theorems as well as
a (more general) $q$-Euler transformation.
Taking the $(q,t)\to (1^{-},1^{-})$ limit of the $\gsl_3$ $q$-binomial theorem
(such that $(1-t)/(1-q)\to \gamma$) yields a generalization of the
Tarasov--Varchenko integral \eqref{integral} involving the Jack polynomial.

\section{Macdonald polynomials}\label{sec2}

\subsection{Preliminaries}\label{sec21}
Let $\la=(\la_1,\la_2,\dots)$ be a partition,
i.e., $\la_1\geq \la_2\geq \dots$ with finitely many $\la_i$
unequal to zero.
The length and weight of $\la$, denoted by
$l(\la)$ and $\abs{\la}$, are the number and sum
of the non-zero $\la_i$ respectively.
As usual we identify two partitions that differ only in their
string of zeros, so that $(6,3,3,1,0,0)$ and $(6,3,1,1)$ represent the
same partition.
When $\abs{\la}=N$ we say that $\la$ is a partition of $N$,
and the unique partition of zero is denoted by $0$.
The multiplicity of the part $i$ in the partition $\la$ is denoted
by $m_i=m_i(\la)$, and occasionally we will write 
$\la=(1^{m_1} 2^{m_2} \dots)$.

We identify a partition with its Ferrers graph,
defined by the set of points in $(i,j)\in \Z^2$ such that
$1\leq j\leq \la_i$, and further make the usual 
identification between Ferrers graphs and (Young) diagrams
by replacing points by squares.

The conjugate $\la'$ of $\la$ is the partition obtained by
reflecting the diagram of $\la$ in the main diagonal,
so that, in particular, $m_i(\la)=\la_i'-\la_{i+1}'$.
The statistic $n(\la)$ is given by
\begin{equation*}
n(\la)=\sum_{i\geq 1} (i-1)\la_i=
\sum_{i\geq 1}\binom{\la_i'}{2}.
\end{equation*}

The dominance partial order on the set of partitions of $N$ is
defined by $\la\geq \mu$ if
$\la_1+\cdots+\la_i\geq \mu_1+\cdots+\mu_i$ for all $i\geq 1$.
If $\la\geq \mu$ and $\la\neq\mu$ then $\la>\mu$.

If $\la$ and $\mu$ are partitions then $\mu\subseteq\la$
if (the diagram of) $\mu$ is contained in (the diagram of)
$\la$, i.e., $\mu_i\leq\la_i$ for all $i\geq 1$.
If $\mu\subseteq\la$ then the skew-diagram $\la-\mu$
denotes the set-theoretic difference between $\la$ and $\mu$,
i.e., those squares of $\la$ not contained in $\mu$.

Let $s=(i,j)$ be a square in the diagram of $\la$. Then
$a(s)$, $a'(s)$, $l(s)$ and $l'(s)$ are the arm-length, arm-colength,
leg-length and leg-colength of $s$, defined by
\begin{align*}
a(s)&=\la_i-j,  & a'(s)&=j-1 \\
l(s)&=\la'_j-i,  & l'(s)&=i-1.
\end{align*}
This may be used to define the generalized hook-length polynomials
\cite[Equation (VI.8.1)]{Macdonald95}
\begin{subequations}\label{cdef}
\begin{align}
c_{\la}(q,t)&=\prod_{s\in\la}\bigl(1-q^{a(s)}t^{l(s)+1}\bigr), \\
c'_{\la}(q,t)&=\prod_{s\in\la}\bigl(1-q^{a(s)+1}t^{l(s)}\bigr),
\end{align}
\end{subequations}
where the products are over all squares of $\la$. We further set
\begin{equation}\label{bdef}
b_{\la}(q,t)=\frac{c_{\la}(q,t)}{c'_{\la}(q,t)}.
\end{equation}
Observe that if $\la$ contains a single part, say $k$, then
\[
c'_{(k)}(q,t)=(q;q)_k.
\]

For $N$ a nonnegative integer the $q$-shifted factorial 
$(b;q)_N$ is defined as $(b;q)_0=1$ and
\begin{equation}\label{qfac}
(b;q)_N=(1-b)(1-bq)\cdots(1-bq^{N-1}).
\end{equation}
We also need the $q$-shifted factorial for negative 
(integer) values of $N$. This may be obtained from the above by
\begin{equation*}
(b;q)_{-N}=\frac{1}{(bq^{-N};q)_N}.
\end{equation*}
This implies in particular that $1/(q;q)_{-N}=0$ for positive $N$.

The definition \eqref{qfac} may be extended to partitions $\la$ by
\begin{equation*}
(b;q,t)_{\la}=\prod_{s\in\la}\bigl(1-b\, q^{a'(s)}t^{-l'(s)}\bigr) \\
=\prod_{i=1}^{l(\la)}(bt^{1-i};q)_{\la_i}.
\end{equation*}
With this notation the polynomials \eqref{cdef} may be recast as
\cite[Proposition 3.2]{Kaneko96}
\begin{subequations}\label{ccp}
\begin{align}\label{c}
c_{\la}(q,t)&=(t^n;q,t)_{\la}
\prod_{1\leq i<j\leq n}\frac{(t^{j-i};q)_{\la_i-\la_j}}
{(t^{j-i+1};q)_{\la_i-\la_j}}, \\
\label{cp}
c'_{\la}(q,t)&=(qt^{n-1};q,t)_{\la}
\prod_{1\leq i<j\leq n}\frac{(qt^{j-i-1};q)_{\la_i-\la_j}}
{(qt^{j-i};q)_{\la_i-\la_j}},
\end{align}
\end{subequations}
where $n$ is any integer such that $n\geq l(\la)$.

Finally we introduce the usual condensed notation for $q$-shifted factorials
as
\begin{equation*}
(a_1,\dots,a_k;q)_N=(a_1;q)_N\cdots (a_k;q)_N
\end{equation*}
and
\begin{equation*}
(a_1,\dots,a_k;q,t)_{\la}=(a_1;q,t)_{\la}\cdots (a_k;q,t)_{\la}.
\end{equation*}

\subsection{Macdonald polynomials}
Let $\Symm_n$ denote the symmetric group, acting on
$x=(x_1,\dots,x_n)$ by permuting the $x_i$,
and let $\Lambda_n=\Z[x_1,\dots,x_n]^{\Symm_n}$ and
$\Lambda$ denote the ring of symmetric polynomials in $n$ independent 
variables and the ring of symmetric functions in countably many variables,
respectively.

For $\la=(\la_1,\dots,\la_n)$ a partition 
of at most $n$ parts the monomial symmetric function $m_{\la}$ 
is defined as
\begin{equation*}
m_{\la}(x)=\sum x^{\alpha},
\end{equation*}
where the sum is over all distinct permutations $\alpha$ of
$\la$, and $x^{\alpha}=x_1^{\alpha_1}\cdots x_n^{\alpha_n}$.
For $l(\la)>n$ we set $m_{\la}(x)=0$.
The monomial symmetric functions $m_{\la}$ for $l(\la)\leq n$
form a $\Z$-basis of $\Lambda_n$.

For $r$ a nonnegative integer the power sums $p_r$ are given by
$p_0=1$ and $p_r=m_{(r)}$ for $r>1$. Hence
\begin{equation}\label{powersums}
p_r(x)=\sum_{i\geq 1} x_i^r.
\end{equation}
More generally the power-sum products are defined as
$p_{\la}(x)=p_{\la_1}(x)\cdots p_{\la_n}(x)$.

Following Macdonald we define the scalar product
$\langle \cdot,\cdot \rangle_{q,t}$ by
\begin{equation*}
\langle p_{\la},p_{\mu}\rangle_{q,t}=
\delta_{\la\mu} z_{\la} \prod_{i=1}^n
\frac{1-q^{\la_i}}{1-t^{\la_i}},
\end{equation*}
with $z_{\la}=\prod_{i\geq 1} m_i! \: i^{m_i}$ and
$m_i=m_i(\la)$. If we denote the ring of symmetric functions 
in $n$ variables over the field $\F=\Q(q,t)$ of rational functions 
in $q$ and $t$ by $\Lambda_{n,\F}$, then
the Macdonald polynomial $P_{\la}(x;q,t)$
is the unique symmetric polynomial in $\Lambda_{n,\F}$ such that
\cite[Equation (VI.4.7)]{Macdonald95}:
\begin{equation}\label{Pm}
P_{\la}(x;q,t)=m_{\la}(x)+\sum_{\mu<\la}
u_{\la\mu}(q,t) m_{\mu}(x)
\end{equation}
and
\begin{equation*}
\langle P_{\la},P_{\mu} \rangle_{q,t}
=0\quad \text{if$\quad\la\neq\mu$.}
\end{equation*}
The Macdonald polynomials $P_{\la}(x;q,t)$ with $l(\la)\leq n$ 
form an $\F$-basis of $\Lambda_{n,\F}$. If $l(\la)>n$ then
$P_{\la}(x;q,t)=0$.
{}From \eqref{Pm} it follows that $P_{\la}(x;q,t)$ for
$l(\la)\leq n$ is homogeneous of degree $\abs{\la}$:
\begin{equation}\label{hom}
P_{\la}(zx;q,t)=z^{\abs{\la}}P_{\la}(x;q,t)
\end{equation}
with $z$ a scalar.

When $q=t$ the Macdonald polynomials simplify to the 
well-known Schur functions:
\begin{equation}\label{qist}
P_{\la}(x;t,t)=s_{\la}(x).
\end{equation}
The latter are defined much more simply as
\begin{equation}\label{Schur}
s_{\la}(x)=
\frac{\det_{1\leq i,j\leq n}\bigl(x_i^{\la_j+n-j}\bigr)}
{\det_{1\leq i,j\leq n}\bigl(x_i^{n-j}\bigr)}=
\frac{\det_{1\leq i,j\leq n}\bigl(x_i^{\la_j+n-j}\bigr)}{\Delta(x)}, 
\end{equation}
where
\begin{equation*}
\Delta(x)=\prod_{1\leq i<j\leq n}(x_i-x_j)
\end{equation*}
is the Vandermonde product.

For $f\in\Lambda_{n,\F}$ and $\la$ a partition such that $l(\la)\leq n$
the evaluation homomorphism 
$u_{\la}^{(n)}:\Lambda_{n,\F}\to \F$ is defined as
\begin{equation}\label{eval}
u_{\la}^{(n)}(f)=f(q^{\la_1}t^{n-1},q^{\la_2}t^{n-2},\dots,q^{\la_n}t^0).
\end{equation}
We extend this to 
$f\in\F(x_1,\dots,x_n)^{\Symm_n}$ for those $f$ for which the
right-hand side of \eqref{eval} is well-defined.
According to the principal specialization formula for Macdonald polynomials 
\cite[Example VI.6.5]{Macdonald95}
\begin{equation}\label{PS}
u_0^{(n)}(P_{\la})
=t^{n(\la)} \prod_{s\in\la}
\frac{1-q^{a'(s)}t^{n-l'(s)}}{1-q^{a(s)}t^{l(s)+1}}
=t^{n(\la)}\frac{(t^n;q,t)_{\la}}{c_{\la}(q,t)}.
\end{equation}
For more general evaluations we have the symmetry
\cite[Equation (VI.6.6)]{Macdonald95}
\begin{equation}\label{symm}
u_{\la}^{(n)}(P_\mu)u_0^{(n)}(P_{\la})=
u_{\mu}^{(n)}(P_{\la})u_0^{(n)}(P_{\mu})
\end{equation}
provided $l(\la),l(\mu)\leq n$. 
It will also be convenient to define the homomorphism $u_{\la;z}^{(n)}$ as
\begin{equation}\label{evalz}
u_{\la;z}^{(n)}(f)=
f(zq^{\la_1}t^{n-1},zq^{\la_2}t^{n-2},\dots,zq^{\la_n}t^0).
\end{equation}
For homogeneous functions of degree $d$ we of course have
\begin{equation}\label{deg}
u_{\la;z}^{(n)}(f)=z^d \, u_{\la}^{(n)}(f).
\end{equation}

Thanks to the stability $P_{\la}(x_1,\dots,x_n;q,t)=
P_{\la}(x_1,\dots,x_n,0;q,t)$ for $l(\la)\leq n$,
we may extend the $P_{\la}$ to an infinite alphabet, and
in the remainder of this section
we assume that $x$ (and $y$) contain countable many variables
so that we will be working in the ring 
$\Lambda_{\F}=\Lambda\otimes_{\Z}\F$ instead of
$\Lambda_{n,\F}$.
By abuse of terminology we still refer to $P_{\la}(x;q,t)$
as a Macdonald polynomial, instead of a Macdonald function.

For $b$ an indeterminate, the homomorphism $\epsilon_{a,t}:\Lambda_{\F}\to \F$
is defined by its action on the power sums $p_r$ 
as \cite[Equation (VI.6.16)]{Macdonald95}
\begin{equation}\label{epsdef}
\epsilon_{b,t}(p_r)=\frac{1-b^r}{1-t^r}.
\end{equation}
According to \cite[Equation (VI.6.17)]{Macdonald95}
\begin{equation}\label{epsP}
\epsilon_{b,t}(P_{\la})
=t^{n(\la)} \prod_{s\in\la}
\frac{1-b\, q^{a'(s)}t^{-l'(s)}}{1-q^{a(s)}t^{l(s)+1}}
=t^{n(\la)}\frac{(b;q,t)_{\la}}{c_{\la}(q,t)}.
\end{equation}
We also note that for any symmetric function $f$
\begin{equation}\label{PS2}
\epsilon_{t^n,t}(f)=u_0^{(n)}(f)=f(1,t,\dots,t^{n-1}),
\end{equation}
compare for example \eqref{PS} and \eqref{epsP}.

The $q,t$-Littlewood--Richardson coefficients
are defined by
\begin{equation}\label{qtLR}
P_{\mu}(x;q,t)P_{\nu}(x;q,t)=
\sum_{\la}f_{\mu\nu}^{\la}(q,t) P_{\la}(x;q,t),
\end{equation}
and trivially satisfy
\begin{equation*}
f_{\mu\nu}^{\la}(q,t)=f_{\nu\mu}^{\la}(q,t)
\end{equation*}
and
\begin{equation}\label{SA}
f_{\mu\nu}^{\la}(q,t)=0 \text{ unless }\abs{\la}=\abs{\mu}+\abs{\nu}.
\end{equation}
It can also be shown that 
\cite[Equation (VI.7.7)]{Macdonald95}
\begin{equation}\label{include}
f_{\mu\nu}^{\la}(q,t)=0 \text{ unless } \mu,\nu\subseteq\la.
\end{equation}

The $q,t$-Littlewood--Richardson coefficients may be used to define
the skew Macdonald polynomials
\begin{equation}\label{skewdef}
P_{\la/\mu}(x;q,t)=\sum_{\nu}f_{\mu\nu}^{\la}(q,t)P_{\nu}(x;q,t).
\end{equation}
By \eqref{include}, $P_{\la/\mu}(x;q,t)=0$
unless $\mu\subseteq\la$ (in which case it is a homogeneous
of degree $\abs{\la}-\abs{\mu}$).
Equivalent to \eqref{skewdef} is
\begin{equation}\label{skewdef2}
P_{\la}(x,y;q,t)=\sum_{\mu} P_{\la/\mu}(x;q,t)P_{\mu}(y;q,t).
\end{equation}

Finally we need the Kaneko--Macdonald definition
of $\gsl_2$ basic hypergeometric series with Macdonald polynomial
argument \cite{Kaneko96,Macdonald}
\begin{equation}\label{Phirs}
{_{r+1}\Phi_r}\hyp{a_1,\dots,a_{r+1}}{b_1,\dots,b_r}{q,t;x}
=\sum_{\lambda} 
t^{n(\la)}\frac{P_{\lambda}(x;q,t)}{c'_{\lambda}(q,t)}\,
\frac{(a_1,\dots,a_{r+1};q,t)_{\lambda}}
{(b_1,\dots,b_r;q,t)_{\lambda}}.
\end{equation}
In the single-variable case, $x=(z)$, 
this reduces to the classical $_{r+1}\phi_r$
basic hypergeometric series \cite{GR04}:
\[
{_{r+1}}\Phi_r\hyp{a_1,\dots,a_{r+1}}{b_1,\dots,b_r}{q,t;(z)}
=\sum_{k=0}^{\infty}
\frac{(a_1,\dots,a_{r+1};q)_k}
{(q,b_1,\dots,b_r;q)_k} \, z^k 
=:{_{r+1}}\phi_r\hyp{a_1,\dots,a_{r+1}}{b_1,\dots,b_r}{q,z}.
\]
The main result for Kaneko--Macdonald series needed in this paper
is the $q$-binomial theorem
\cite[Theorem 3.5]{Kaneko96}, \cite[Equation (2.2)]{Macdonald}
(see also \cite[page 374]{Macdonald95})
\begin{equation}\label{Phi10sum}
{_1\Phi_0}\hyp{a}{\text{--}}{q,t;x}
=\prod_{i\geq 1} \frac{(ax_i;q)_{\infty}}{(x_i;q)_{\infty}}
\end{equation}
which is \eqref{qbthm2} of the introduction.
Those familiar with Macdonald polynomials will recognize the
intimate connection with the
Cauchy identity \cite[Equation (VI.4.13)]{Macdonald95}
\begin{equation}\label{Cauchy}
\sum_{\la} b_{\la}(q,t)P_{\la}(x;q,t)P_{\la}(y;q,t)
=\prod_{i,j\geq 1} \frac{(tx_i y_j;q)_{\infty}}{(x_i y_j;q)_{\infty}},
\end{equation}
with $b_{\la}(q,t)$ defined in \eqref{bdef}.
Acting with the homomorphism $\epsilon_{a,t}$ on the left
(with $\epsilon_{a,t}$ acting on $y$)
and using \eqref{bdef} and \eqref{epsP} immediately gives the
above $_1\Phi_0$ series, so that \eqref{Phi10sum} is equivalent
to 
\begin{equation}\label{epsqbin}
\epsilon_{a,t}\biggl(\:
\prod_{i,j\geq 1} 
\frac{(tx_i y_j;q)_{\infty}}{(x_i y_j;q)_{\infty}} \biggr)
=\prod_{i\geq 1} \frac{(ax_i;q)_{\infty}}{(x_i;q)_{\infty}}.
\end{equation}

\section{The bisymmetric function $F$}\label{secbiF}
Unless stated otherwise $m$ and $n$ are integers such that
$0\leq m\leq n$, and $x=(x_1,\dots,x_m)$ and $y=(y_1,\dots,y_n)$.
Given such $x$ we set 
\[
x^{(i_1,i_2,\dots,i_N)}=(x_1,\dots,x_{i_1-1},x_{i_1+1},\dots,x_{i_2-1},
x_{i_2+1},\dots,x_{i_N-1},x_{i_N+1},\dots,x_m)
\]
for integers $1\leq i_1<i_2<\dots<i_N\leq m$.
We further use the shorthand notation 
\begin{equation*}
(x^{(p+1,\dots,m)},0^{m-p})=(x_1,\dots,x_p,
\underbrace{0,\dots,0}_{m-p~\text{times}}\!\!),
\end{equation*}
and apply the same notation to $y=(y_1,\dots,y_n)$.

The symmetric group will feature prominently in this section,
especially in the proofs. In total we employ the symmetric group
acting on $4$ different sets of variables, sometimes of the same
cardinality. To avoid ambiguity we write 
\[
\sum_{w\in\Symm_x}w\bigl(f(x)\bigr)
\]
instead of the more common
\[
\sum_{w\in\Symm_m}w\bigl(f(x)\bigr):=
\sum_{w\in\Symm_m}f(x_{w_1},\dots,x_{w_m}),
\]
with similar notation for other sets of variables.

\subsection{Definitions and results}\label{secDR}

Let $r$ be a nonnegative integer not exceeding $m$.
Macdonald introduced the commuting family of $q$-difference operators 
$D_r$ as \cite[Equation (VI.3.4)$_r$]{Macdonald95}
\[
D_r=t^{\binom{r}{2}}\sum_{\substack{I\subseteq [m] \\ \abs{I}=r}}
\prod_{\substack{i\in I \\ j\not\in I}}
\frac{tx_i-x_j}{x_i-x_j} \prod_{i\in I} T_{q,x_i},
\]
where $[m]=\{1,2,\dots,m\}$ and
\[
T_{q,x_i}\bigl(f(x)\bigr)
=f(x_1,\dots,x_{i-1},qx_i,x_{i+1},\dots,x_m)
\]
the $q$-shift operator acting on $x_i$.

Defining the generating series 
\[
D(u;q,t)=\sum_{r=0}^m D_r (-u)^r
\]
Macdonald showed that for $l(\la)\leq m$ the $P_{\la}$ are
the eigenfunctions of $D(u;q,t)$
\cite[Equation (VI.4.15)]{Macdonald95}:
\begin{equation}\label{Lambda}
D(u;q,t) P_{\la}(x;q,t)=g_{\la}(u;q,t) P_{\la}(x;q,t),
\end{equation}
with eigenvalue
\[
g_{\la}(u;q,t)=\prod_{i=1}^m (1-ut^{m-i}q^{\la_i}).
\]

In \cite[Equations (1.12) and (1.13)]{KN99} Kirillov and Noumi combined 
the Cauchy identity \eqref{Cauchy} with \eqref{Lambda} to obtain 
\begin{equation}\label{Fdef}
\sum_{\la} b_{\la}(q,t) g_{\la}(u;q,t)
P_{\la}(x;q,t)P_{\la}(y;q,t)
=F(u;x,y;t) \prod_{i=1}^m \prod_{j=1}^n 
\frac{(tx_i y_j;q)_{\infty}}{(x_i y_j;q)_{\infty}},
\end{equation}
where the bisymmetric function $F(u;x,y;t)$ is given by
\begin{equation}\label{Fexp}
F(u;x,y;t)=\sum_{I\subseteq [m]}
(-u)^{\abs{I}}t^{\binom{\abs{I}}{2}}
\prod_{\substack{i\in I \\ j\not\in I}}
\frac{tx_i-x_j}{x_i-x_j}
\prod_{i\in I}\prod_{j=1}^n
\frac{1-x_iy_j}{1-tx_iy_j}.
\end{equation}

In the Section~\ref{SecA2} we define two types of $\gsl_3$
basic hypergeometric series featuring particular
specializations of $F$. In our study of these series several
elementary results for $F$ are needed. Proofs of all claims
may be found in Section~\ref{secFProofs}.

\begin{lemma}[(Stability)]\label{lemS1}
We have
\begin{subequations}
\begin{align}\label{stab1}
F(u;x,y;t)|_{x_my_n=1}\:&=\:F(u;x^{(m)},y^{(n)};t) \\
\intertext{and}
F(u;x,y;t)|_{x_m=y_n=0}\:&=\:(1-u)F(ut;x^{(m)},y^{(n)};t).
\label{stab2}
\end{align}
\end{subequations}
\end{lemma}

The formulae \eqref{Fdef} and \eqref{Fexp} also make sense when
$y$ contains countably many variables
(provided, of course, we replace $\prod_{j=1}^n$ by $\prod_{j\geq 1}$).
In the following we assume such $y$.
\begin{lemma}\label{lemeps}
With $\epsilon_{a,t}$ acting on $y=(y_1,y_2,\dots)$ we have
\begin{subequations}
\begin{equation}\label{E1}
\epsilon_{ut^{m-1},t}\bigl(F(u;x,y;t)\bigr)
=\prod_{i=1}^m \frac{1-ut^{m-i}}{1-ut^{m-1}x_i}
\end{equation}
and
\begin{equation}\label{E2}
\epsilon_{a,t}\bigl(F(1;x,y;t)\bigr) \\
=t^{\binom{m}{2}}x_1\cdots x_m\prod_{i=1}^m \frac{1-at^{1-i}}{1-ax_i}.
\end{equation}
\end{subequations}
\end{lemma}
It easily follows (see Section~\ref{secFProofs}) that
\begin{equation}\label{epsF}
\epsilon_{a,t}\bigl(F(u;x,y;t)\bigr)=\sum_{I\subseteq [m]}
(-u)^{\abs{I}}t^{\binom{\abs{I}}{2}}
\prod_{\substack{i\in I \\ j\not\in I}}
\frac{tx_i-x_j}{x_i-x_j}
\prod_{i\in I} \frac{1-x_i}{1-ax_i},
\end{equation}
so that Lemma~\ref{lemeps} is equivalent to the pair of identities
\begin{subequations}
\begin{equation}\label{E1b}
\sum_{I\subseteq [m]}
(-u)^{\abs{I}}t^{\binom{\abs{I}}{2}}
\prod_{\substack{i\in I \\ j\not\in I}}
\frac{tx_i-x_j}{x_i-x_j}
\prod_{i\in I} \frac{1-x_i}{1-ut^{m-1}x_i}
=\prod_{i=1}^m \frac{1-ut^{m-i}}{1-ut^{m-1}x_i}
\end{equation}
and
\begin{equation}\label{E2b}
\sum_{I\subseteq [m]}
(-1)^{\abs{I}}t^{\binom{\abs{I}}{2}}
\prod_{\substack{i\in I \\ j\not\in I}}
\frac{tx_i-x_j}{x_i-x_j}
\prod_{i\in I} \frac{1-x_i}{1-ax_i}
=t^{\binom{m}{2}}x_1\cdots x_m\prod_{i=1}^m \frac{1-at^{1-i}}{1-ax_i}.
\end{equation}
\end{subequations}
This shows that \eqref{E1} and \eqref{E2} are in fact equivalent:
taking \eqref{E1b} and making the substitutions
$u\to at^{m-1}$, $x_i\to 1/(a x_i)$ and $I\to [m]-I$ yields \eqref{E2b}.

The results that we will actually need in Section~\ref{SecA2} 
correspond to the principal specialization formula, obtained 
by choosing $u=t^{n-m+1}$ or $a=t^n$ in
Lemma~\ref{lemeps} and using \eqref{PS2}.

\begin{corollary}[(Principal specialization)]\label{corPSF}
With $u_0^{(n)}$ acting on $y=(y_1,\dots,y_n)$ we have
\[
u_0^{(n)}\bigl(F(t^{n-m+1};x,y;t)\bigr)=
\prod_{i=1}^m \frac{1-t^{i+n-m}}{1-t^n x_i}
\]
and
\[
u_0^{(n)}\bigl(F(1;x,y;t)\bigr) \\
=t^{\binom{m}{2}}x_1\cdots x_m
\prod_{i=1}^m \frac{1-t^{i+n-m}}{1-t^n x_i}.
\]
\end{corollary}
These last two results are suggestive of
\[
F(1;x,y;t)=t^{\binom{m}{2}-\binom{n}{2}}
F(t^{n-m+1};x,y;t)\prod_{i=1}^m x_i\prod_{j=1}^n y_j,
\]
but this is in fact only true for $m=n$ as will be shown
in \eqref{F1t} below.

\medskip

The function $F$ may be connected to the bisymmetric function
introduced by Tarasov and Varchenko in their work on $\gsl_3$
Selberg integrals \cite{TV03}.
To this end we define 
\begin{equation}\label{omegaF}
\omega(x,y;t)=F(1;x^{-1},y;t),
\end{equation}
where $x^{-1}=(x_1^{-1},\dots,x_m^{-1})$.
{}From \eqref{Fexp} it follows that
\begin{equation}\label{omega}
\omega(x,y;t)=\sum_{I\subseteq [m]}
(-1)^{\abs{I}}t^{\binom{\abs{I}}{2}}
\prod_{\substack{i\in I \\ j\not\in I}} \frac{x_i-tx_j}{x_i-x_j}
\prod_{i\in I}\prod_{j=1}^n \frac{x_i-y_j}{x_i-ty_j}.
\end{equation}

\begin{proposition}\label{proprec} 
Let $k$ be an integer such that $1\leq k\leq m$. Then
\begin{equation}\label{it}
\omega(x,y;t)=t^{m-n}(1-t)\sum_{l=1}^n\omega(x^{(k)},y^{(l)};t)\,
\frac{y_l}{x_k-ty_l}  
\prod_{\substack{i=1\\ i\neq k}}^m \frac{x_i-y_l}{x_i-ty_l} 
\prod_{\substack{i=1 \\ i\neq l}}^n \frac{y_i-ty_l}{y_i-y_l}.
\end{equation}
\end{proposition}

Since $\omega(\text{--}\,,y;t)=1$ we may use \eqref{it}
and induction to find the following alternative
multisum expression for $\omega$.
\begin{corollary}\label{corom}
We have
\begin{multline}\label{itt}
\omega(x,y;t)=t^{m(m-n)}(1-t)^m \\ \times
\sum_{\substack{l_1,\dots,l_m=1 \\ l_i\neq l_j}}^n
\prod_{i=1}^m \frac{y_{l_i}}{x_i-ty_{l_i}}
\prod_{\substack{i=1 \\ i\neq l_1,\dots,l_m}}^n
\prod_{j=1}^m\frac{y_i-ty_{l_j}}{y_i-y_{l_j}} 
\prod_{1\leq i<j\leq m}
\frac{x_i-y_{l_j}}{x_i-ty_{l_j}}\cdot
\frac{y_{l_i}-ty_{l_j}}{y_{l_i}-y_{l_j}}.
\end{multline}
\end{corollary}
Note that for $m=n$ this is equivalent to
\[
\omega(x,y;t)=(1-t)^n
\sum_{w\in\Symm_y}w\biggl(\;
\prod_{i=1}^n \frac{y_i}{x_i-ty_i}
\prod_{1\leq i<j\leq n}
\frac{x_i-y_j}{x_i-ty_j}\cdot
\frac{y_i-ty_j}{y_i-y_j}\biggr)
\]
from which it readily follows that
\[
\omega(x,y;t)=\omega(x^{-1},y^{-1};t^{-1})\prod_{i=1}^n\frac{y_i}{x_i}
\]
or, equivalently,
\[
F(1;x,y;t)=F(1;x^{-1},y^{-1};t^{-1})\prod_{i=1}^n x_iy_i.
\]
Since it follows from \eqref{Fexp} that for general 
$0\leq m\leq n$
\[
F(u;x,y;t)=F(ut^{m-n-1};x^{-1},y^{-1};t^{-1}),
\]
we also have
\begin{equation}\label{F1t}
F(1;x,y;t)=F(t;x,y;t)\prod_{i=1}^n x_iy_i
\end{equation}
when $m=n$.

Using Corollary~\ref{corom} we may achieve the further rewriting
of $\omega$ as follows.
\begin{proposition}\label{propomegaSS}
We have
\begin{multline}\label{omegaSS}
\omega(x,y;t)=
\frac{t^{m(m-n)} (1-t)^{n+m}}{(t;t)_{n-m}(t;t)_m} 
\sum_{w\in \Symm_x\times\Symm_y}w\biggl(\;
\prod_{i=1}^m\frac{y_{i+n-m}}{x_i-ty_{i+n-m}} 
\prod_{1\leq i<j\leq n} \frac{y_i-ty_j}{y_i-y_j} \\
\times \prod_{1\leq i<j\leq m}
\frac{x_i-y_{j+n-m}}{x_i-ty_{j+n-m}}\cdot
\frac{x_i-tx_j}{x_i-x_j} \biggr).
\end{multline}
\end{proposition}
The representation of $\omega(x,y;t)$ provided by \eqref{omegaSS}
immediately implies that
\begin{equation}\label{qdeformed}
\lim_{q\to 1} F(1;q^{-v},q^u;q^{\gamma})=
\lim_{q\to 1} \omega(q^v,q^u;q^{\gamma})=
\frac{(-\gamma)^m n!}{(n-m)!}\, w(u,v;\gamma),
\end{equation}
where $w(u,v;\gamma)$ is the bisymmetric function of
Tarasov and Varchenko \cite[Eq. (2.2)]{TV03}, and
$q^v=(q^{v_1},\dots,q^{v_m})$,
$q^u=(q^{u_1},\dots,q^{u_n})$.

Depending on the respective values of $m$ and $n$
either \eqref{omega} or \eqref{itt}
provides the most efficient way of computing $\omega(x,y;t)$.
In the former we need to sum over all $2^m$ subsets 
of $[m]$ whereas in the latter we are summing over
all $\binom{n}{m}$ $m$-subsets of $[n]$.
A distinct advantage of the representation \eqref{itt}
(and of \eqref{omegaSS}) over \eqref{omega} is that it permits the
computation of the $t\to 1$ limit, required in the derivation of the
$\gsl_3$ Selberg integral \eqref{integral}.
In particular, the bisymmetric function featured in that
integral follows as
\begin{align}\label{gdef}
h(x,y)&=(-1)^m \, \frac{(n-m)!}{n!}
\lim_{t\to 1}\frac{\omega(x,y;t)}{(1-t)^m} \\
&=\frac{(n-m)!}{n!}\sum_{\substack{l_1,\dots,l_m=1 \\ l_i\neq l_j}}^n
\prod_{i=1}^m \frac{y_{l_i}}{y_{l_i}-x_i} \notag  \\
&=\frac{1}{m!n!} 
\sum_{w\in \Symm_x\times\Symm_y}w\biggl(\;
\prod_{i=1}^m\frac{y_{i+n-m}}{y_{i+n-m}-x_i} \biggr).
\notag 
\end{align}

Finally we mention that $F(t;x,y;t)$ for $m=n$ is nothing but the 
well-known Izergin--Korepin determinant \cite{Izergin87,Korepin82} 
in disguise.
\begin{lemma}\label{lemdet}
For $x=(x_1,\dots,x_n)$ and $y=(y_1,\dots,y_n)$ we have
\begin{equation*}
F(t;x,y;t)=
\det_{1\leq i,j\leq n}\biggl(
\frac{1}{(1-x_iy_j)(1-tx_iy_j)}\biggr)
\frac{(1-t)^n\prod_{i,j=1}^n (1-x_i y_j)}
{\prod_{1\leq i<j\leq n}(x_i-x_j)(y_i-y_j)}.
\end{equation*}
\end{lemma}
Since $F(0;x,y;0)=1$ this reduces to
Cauchy's double alternant when $t=0$, see e.g., 
\cite[Equation 2.7]{Krattenthaler99}.

Several combinatorial interpretations of the 
Izergin--Korepin determinant are known, for example as
the partition function of square ice \cite{Bressoud99,Lascoux99}.
Perhaps best known is its evaluation in terms of alternating sign
matrices \cite{Bressoud99,Kuperberg96}. This (together with
\eqref{omegaF} and \eqref{F1t}) implies that for $m=n$ 
\[
\omega(x,y;t)=\frac{(1-t)^n y_1\cdots y_n}
{\prod_{i,j=1}^n(x_i-ty_j)} 
\sum_A(1-t)^{2N(A)}t^{\binom{n}{2}-\mathcal{I}(A)}
\prod_{i=1}^ny_i^{N_i(A)}x_i^{N^i(A)}
\prod_{\substack{i,j=1 \\ a_{ij}=0}}^n
(\alpha_{ij}y_i-x_j).
\]
Here the sum is over all $n$ by $n$ alternating sign matrices $A$
(matrices with entries $a_{ij}\in\{-1,0,1\}$ such that the 
ones and minus ones
alternate along each row and along each column and such that the entries in
each row and column add up to $1$), $N_i(A)$ is the number of
minus ones in row $i$, $N^i(A)$ is the number of minus ones in column $i$,
$N(A)$ is the total number of minus ones,
$\mathcal{I}(A)$ is the inversion number:
\begin{equation*}
\mathcal{I}(A)=
\sum_{1\leq i'<i\leq n}\sum_{1\leq j<j'\leq n}a_{ij}a_{i'j'},
\end{equation*}
and 
\begin{equation*}
\alpha_{ij}=t\quad \text{ if }\quad \sum_{k=1}^ja_{ik}=\sum_{k=1}^ia_{kj}
\end{equation*}
and $\alpha_{ij}=1$ otherwise.

\subsection{The rational functions $W_{\la\mu}$ and $V_{\la\mu}$}

Related to the bisymmetric function $F$ we introduce
two rational functions $W_{\la\mu}(u,z;q,t)$
and $V_{\la\mu}(u,z;q,t)$ as follows.
Let $\la$ and $\mu$ be partitions such that $l(\la)\leq m$ and
$l(\mu)\leq n$. Then
\begin{equation}\label{Wlamudef}
W_{\la\mu}(u,z;q,t)=u_{\la;z}^{(m)}u_{\mu}^{(n)}\bigl(F(u;x,y;t)\bigr)
\end{equation}
and
\begin{equation}\label{Vlamudef}
V_{\la\mu}(u,z;q,t)=
u_{\la;z}^{(m)}u_{\mu}^{(n)}\bigl(F(u;x^{-1},y;t)\bigr).
\end{equation}
There is no need to consider the more general specialization
$\displaystyle u_{\la;z}^{(m)}u_{\mu;w}^{(n)}$
since
\[
u_{\la;z}^{(m)}u_{\mu;w}^{(n)}\bigl(F(u;x,y;t)\bigr)
=u_{\la;zw}^{(m)}u_{\mu}^{(n)}\bigl(F(u;x,y;t)\bigr) .
\]

{}From \eqref{Fexp} it immediately follows that
\[
W_{\la\mu}(u,z;q,t)=\sum_{I\subseteq [m]}
(-u)^{\abs{I}}t^{\binom{\abs{I}}{2}}
\prod_{\substack{i\in I \\ j\not\in I}}
\frac{1-q^{\la_i-\la_j}t^{j-i+1}}{1-q^{\la_i-\la_j}t^{j-i}} 
\prod_{i\in I}\prod_{j=1}^n
\frac{1-zq^{\la_i+\mu_j}t^{m+n-i-j}}
{1-zq^{\la_i+\mu_j}t^{m+n-i-j+1}}
\]
and
\[
V_{\la\mu}(u,z;q,t)=\sum_{I\subseteq [m]}
(-u)^{\abs{I}} t^{\binom{\abs{I}}{2}-n\abs{I}}
\prod_{\substack{i\in I \\ j\not\in I}}
\frac{1-q^{\la_j-\la_i}t^{i-j+1}}{1-q^{\la_j-\la_i}t^{i-j}} 
\prod_{i\in I}\prod_{j=1}^n
\frac{1-zq^{\la_i-\mu_j}t^{j-i+m-n}}{1-zq^{\la_i-\mu_j}t^{j-i+m-n-1}}.
\]
{}Furthermore, from 
\eqref{Vlamudef} and Corollary~\ref{corPSF} we infer that
\begin{subequations}
\begin{equation}\label{VPS}
V_{\la,0}(t^{n-m+1},z;q,t)=
q^{\abs{\la}} z^m 
\prod_{i=1}^m \frac{1-t^{m-n-i}}{1-zq^{\la_i}t^{m-n-i}}
\end{equation}
and
\begin{equation}\label{VPS2}
V_{\la,0}(1,z;q,t)=
\prod_{i=1}^m \frac{1-t^{m-n-i}}{1-zq^{\la_i}t^{m-n-i}}.
\end{equation}
\end{subequations}

\subsection{Proofs of the claims of Section~\ref{secDR}}\label{secFProofs}

\begin{proof}[Proof of Lemma~\ref{lemS1}]
By taking $x_m y_n=1$ in \eqref{Fexp} it follows
that the summand vanishes if $m\in I$.
Hence we need to only sum over $I\subseteq[m-1]$, resulting in
\begin{align*}
F(u;x,y;t)|_{x_m y_n=1} 
&=\sum_{I\subseteq [m-1]}
(-u)^{\abs{I}}t^{\binom{\abs{I}}{2}}
\prod_{\substack{i\in I \\ j\in [m]-I}}
\frac{tx_i-x_j}{x_i-x_j} 
\prod_{i\in I}\biggl(
\frac{1-x_i/x_m}{1-tx_i/x_m}
\prod_{j=1}^{n-1}
\frac{1-x_iy_j}{1-tx_iy_j}\biggr) \\
&=\sum_{I\subseteq [m-1]}
(-u)^{\abs{I}}t^{\binom{\abs{I}}{2}}
\prod_{\substack{i\in I \\ j\not\in I}}
\frac{tx_i-x_j}{x_i-x_j}
\prod_{i\in I} \prod_{j=1}^{n-1}
\frac{1-x_iy_j}{1-tx_iy_j}.
\end{align*}
This last expression is $F(u;x^{(m)},y^{(n)};t)$, establishing 
\eqref{stab1}.

In proving \eqref{stab2} we make the $m$-dependence of
$g_{\la}(u;q,t)$ explicit by writing $g^{(m)}_{\la}(u;q,t)$.

Taking $x_m=y_n=0$ in \eqref{Fdef} and using the stability of the
Macdonald polynomials yields
\[
\sum_{\la} b_{\la}(q,t) g^{(m)}_{\la}(u;q,t)
P_{\la}(x^{(m)};q,t)P_{\la}(y^{(n)};q,t) 
=F(u;x,y;t)|_{x_m=y_n=0} \prod_{i=1}^{m-1} \prod_{j=1}^{n-1} 
\frac{(tx_i y_j;q)_{\infty}}{(x_i y_j;q)_{\infty}}.
\]
Since $P_{\la}(x^{(m)};q,t)=0$ if $l(\la)\geq m$
we may assume that $l(\la)\leq m-1$. But then
\begin{align*}
g^{(m)}_{\la}(u;q,t)&=(1-u)\prod_{i=1}^{m-1} (1-ut^{m-i}q^{\la_i}) \\
&=(1-u)\, g^{(m-1)}_{\la}(ut;q,t),
\end{align*}
so that
\begin{multline*}
(1-u)\sum_{\la} b_{\la}(q,t) g^{(m-1)}_{\la}(ut;q,t)
P_{\la}(x^{(m)};q,t)P_{\la}(y^{(n)};q,t) \\
=F(u;x,y;t)|_{x_m=y_n=0} \prod_{i=1}^{m-1} \prod_{j=1}^{n-1} 
\frac{(tx_i y_j;q)_{\infty}}{(x_i y_j;q)_{\infty}}.
\end{multline*}
Summing the left-hand side using \eqref{Fdef} (with $(n,m,x,y)\to
(n-1,m-1,x^{(m)},y^{(n)})$) completes the proof of \eqref{stab2}.
\end{proof}

\begin{proof}[Proof of Lemma~\ref{lemeps}]

Recall our earlier comment following \eqref{Phi10sum} that the
$_1\Phi_0$ series naturally arises from the sum side of
the Cauchy identity \eqref{Cauchy} by application
of the homomorphism $\epsilon_{a,t}$ (acting on $y$).
It is therefore an obvious idea to apply 
$\epsilon_{a,t}$ to the more more general identity 
\[
\sum_{\la} b_{\la}(q,t) g_{\la}(u;q,t)
P_{\la}(x;q,t)P_{\la}(y;q,t)
=F(u;x,y;t) \prod_{i=1}^m \prod_{j=1}^{\infty}
\frac{(tx_i y_j;q)_{\infty}}{(x_i y_j;q)_{\infty}}.
\]
Doing so and using \eqref{bdef}, \eqref{epsP}, \eqref{Phirs},
\eqref{epsqbin} and 
\[
g_{\la}(u;q,t)=g_0(u;q,t)\,
\frac{(uqt^{m-1};q,t)_{\la}}{(ut^{m-1};q,t)_{\la}},
\]
yields
\begin{equation}\label{phi21}
g_0(u;q,t)\;
{_2\Phi_1}\hyp{a,uqt^{m-1}}{ut^{m-1}}{q,t;x}
=
\epsilon_{a,t}\bigl(F(u;x,y;t)\bigr) 
\prod_{i=1}^m \frac{(ax_i;q)_{\infty}}{(x_i;q)_{\infty}}
\end{equation}
or, equivalently,
\begin{equation}\label{agen}
\epsilon_{a,t}\bigl(F(u;x,y;t)\bigr) \\
= g_0(u;q,t)\;
{_2\Phi_1}\hyp{a,uqt^{m-1}}{ut^{m-1}}{q,t;x}
\prod_{i=1}^m \frac{(x_i;q)_{\infty}}{(ax_i;q)_{\infty}}.
\end{equation}
Taking $a=ut^{m-1}$ the ${_2\Phi_1}$ reduces to a
${_1\Phi_0}$ which may be summed by \eqref{Phi10sum}, so that
\begin{align*}
\epsilon_{ut^{m-1},t}\bigl(F(u;x,y;t)\bigr) &=
g_0(u;q,t)
\prod_{i=1}^m \frac{(uqt^{m-1}x_i;q)_{\infty}}{(ut^{m-1}x_i;q)_{\infty}}\\
&=\prod_{i=1}^m \frac{1-ut^{m-i}}{1-ut^{m-1}x_i}
\end{align*}
in accordance with \eqref{E1}.

To also prove \eqref{E2} we have to prove 
identity \eqref{epsF} (see the comments
immediately following Lemma~\ref{lemeps}).
Hence we need to show that
\[
\epsilon_{a,t}\biggl(\;\prod_{j\geq 1} \frac{1-zy_j}{1-tzy_j}\biggr)=
\frac{1-z}{1-az}.
\]
By taking the logarithm on both sides this is equivalent to
\[
\epsilon_{a,t}\biggl(\;\sum_{j\geq 1}
\Bigl(\log(1-zy_j)-\log(1-tzy_j)\Bigr)\biggr)=
\log\Bigl(\frac{1-z}{1-az}\Bigr).
\]
Using the series expansion for $\log(1-x)$, then
interchanging sums and finally
using definition \eqref{powersums} of the power sums,
this yields
\[
\epsilon_{a,t}\biggl(-\sum_{m\geq 1}\frac{(1-t^m) z^m}{m}\: p_m(y)
\biggr)=\log\Bigl(\frac{1-z}{1-az}\Bigr).
\]
By \eqref{epsdef} this simplifies to
\[
-\sum_{m\geq 1}\frac{(1-a^m) z^m}{m}=\log\Bigl(\frac{1-z}{1-az}\Bigr)
\]
which is obviously true.
\end{proof}

As an aside we note that \eqref{epsF} and \eqref{phi21} 
may be combined to yield
the following generalization of the Kaneko--Macdonald
$q$-binomial theorem \eqref{Phi10sum}:
\begin{multline*}
{_2\Phi_1}\hyp{a,uqt^{m-1}}{ut^{m-1}}{q,t;x}
\prod_{i=1}^m (1-u t^{m-i}) \\
=\biggl(\:\prod_{i=1}^m 
\frac{(ax_i;q)_{\infty}}{(x_i;q)_{\infty}} \biggr)
\sum_{I\subseteq [m]}
(-u)^{\abs{I}}t^{\binom{\abs{I}}{2}}
\prod_{\substack{i\in I \\ j\in \bar{I}}}
\frac{tx_i-x_j}{x_i-x_j}
\prod_{i\in I} \frac{1-x_i}{1-ax_i}.
\end{multline*}

\begin{proof}[Proof of Proposition~\ref{proprec}]
Since $\omega(x,y;t)$ is symmetric in $x$ it suffices to
prove the proposition for $k=m$.

It follows from \eqref{omega} that $\omega(x,y;t)$, viewed as a 
function of $x_m$, has simple poles at $x_m=x_i$ for $1\leq i\leq m-1$
and $x_m=ty_j$ for $1\leq j\leq n$.
However, since $\omega(x,y;t)$ is symmetric in $x$,
the first set of poles must have zero residue. 

It also follows from \eqref{omega} that
\[
\lim_{x_m\to\infty} \omega(x,y;t)=0.
\]
Indeed, if $\omega_I(x,y;t)$ is the summand of \eqref{omega}
and if $I\subseteq [m-1]$, then
\[\lim_{x_m\to\infty} \omega_I(x,y;t)=
-\lim_{x_m\to\infty} \omega_{I\cup\{m\}}(x,y;t).
\]
The above observations imply the existence of the partial 
fraction expansion
\[
\omega(x,y;t)=\sum_{l=1}^n \frac{A_l}{x_m-ty_l},
\]
with $A_l=A_l(x^{(m)},y;t)$ determined by
\begin{align*}
A_l&=\lim_{x_m\to t y_l}
(x_m-ty_l)\, \omega(x,y;t) \\
&=\lim_{x_m\to t y_l}(x_m-ty_l)
\sum_{I\subseteq [m]}
(-1)^{\abs{I}}t^{\binom{\abs{I}}{2}}
\prod_{\substack{i\in I \\ j\not\in I}} \frac{x_i-tx_j}{x_i-x_j}
\prod_{i\in I}\prod_{j=1}^n \frac{x_i-y_j}{x_i-ty_j}.
\end{align*}
In the limit, only sets $I$ containing $m$ give a nonvanishing contribution.
A straightforward calculation thus gives
\begin{multline*}
A_l=(t-1)y_l
\sum_{\substack{I\subseteq [m] \\ m\in I}}
(-1)^{\abs{I}}t^{\binom{\abs{I}}{2}-\abs{I}+m-n+1}
\prod_{j\not\in I} \frac{x_j-y_l}{x_j-ty_l}
\prod_{\substack{j=1 \\ j\neq l}}^n \frac{y_j-ty_l}{y_j-y_l} \\
\times
\prod_{\substack{i\in I-\{m\} \\ j\not\in I}} 
\frac{x_i-tx_j}{x_i-x_j}
\prod_{i\in I-\{m\}}
\prod_{j=1}^n \frac{x_i-y_j}{x_i-ty_j}.
\end{multline*}

Rewriting the sum as a sum over $[m-1]$ this becomes
\[
A_l=(t-1)y_l
\sum_{I\subseteq [m-1]}
(-1)^{\abs{I}+1}t^{\binom{\abs{I}}{2}+m-n}
\prod_{j\not\in I} \frac{x_j-y_l}{x_j-ty_l}
\prod_{\substack{j=1 \\ j\neq l}}^n \frac{y_j-ty_l}{y_j-y_l} 
\prod_{\substack{i\in I \\ j\not\in I}} 
\frac{x_i-tx_j}{x_i-x_j}
\prod_{i\in I} \prod_{j=1}^n \frac{x_i-y_j}{x_i-ty_j}.
\]
By
\[
\prod_{j\not\in I} \frac{x_j-y_l}{x_j-ty_l}
\prod_{i\in I} \prod_{j=1}^n \frac{x_i-y_j}{x_i-ty_j}
=\prod_{i=1}^{m-1} \frac{x_i-y_l}{x_i-ty_l}
\prod_{i\in I} \prod_{\substack{j=1 \\ j\neq l}}^n 
\frac{x_i-y_j}{x_i-ty_j},
\]
this finally yields
\begin{align*}
A_l&=(1-t)t^{m-n}y_l \prod_{i=1}^{m-1} \frac{x_i-y_l}{x_i-ty_l}
\prod_{\substack{j=1 \\ j\neq l}}^n \frac{y_j-ty_l}{y_j-y_l} 
\sum_{I\subseteq [m-1]}
(-1)^{\abs{I}}t^{\binom{\abs{I}}{2}}
\prod_{\substack{i\in I \\ j\in \bar{I}}} 
\frac{x_i-tx_j}{x_i-x_j}
\prod_{i\in I} \prod_{\substack{j=1 \\ j\neq l}}^n 
\frac{x_i-y_j}{x_i-ty_j} \\
&=(1-t)t^{m-n}y_l \,
\omega(x^{(m)},y^{(l)};t)
\prod_{i=1}^{m-1} \frac{x_i-y_l}{x_i-ty_l}
\prod_{\substack{j=1 \\ j\neq l}}^n \frac{y_j-ty_l}{y_j-y_l} 
\end{align*}
as required.
\end{proof}

\begin{proof}[Proof of Proposition~\ref{propomegaSS}]
We first symmetrize the right-hand side of \eqref{omegaSS}
with respect to $y$ and compute
\[
\sum_{w\in \Symm_y} w\biggl(\;
\prod_{i=1}^m\frac{y_{i+n-m}}{x_i-ty_{i+n-m}} 
\prod_{1\leq i<j\leq n}
\frac{y_i-ty_j}{y_i-y_j}
\prod_{1\leq i<j\leq m}
\frac{x_i-y_{j+n-m}}{x_i-ty_{j+n-m}}\biggr).
\]
To this end we write each permutation $w$ 
as $w=(\sigma_1,\dots,\sigma_{n-m},l_1,\dots,l_m)$.
In summing over $w$ we first sum over the $\sigma_i$ for
fixed $l_1,\dots,l_m$. This yields,
\[
\sum_{\substack{l_1,\dots,l_m=1\\ l_i\neq l_j}}^n
\prod_{i=1}^m\frac{y_{l_i}}{x_i-ty_{l_i}} 
\prod_{i=1}^{n-m} 
\prod_{j=1}^m \frac{Y_i-ty_{l_j}}{Y_i-y_{l_j}} 
\prod_{1\leq i<j\leq m}
\frac{x_i-y_{l_j}}{x_i-ty_{l_j}} \cdot 
\frac{y_{l_i}-ty_{l_i}}{y_{l_i}-y_{l_j}}  
\sum_{\sigma\in\Symm_Y} \sigma \biggl(\;
\prod_{1\leq i<j\leq n-m} \frac{Y_i-tY_j}{Y_i-Y_j}
\biggr),
\]
where $Y=(Y_1,\dots,Y_{n-m})=y^{(l_1,l_2,\dots,l_m)}$
and where we have used the symmetry of the double product 
involving $Y_i$ and $y_{l_j}$
to pull it out of the sum over $\Symm_Y$.
Carrying out this sum using \cite[Ch. III, (1.4)]{Macdonald95}
\begin{equation}\label{M14}
\sum_{w\in\Symm_u} w\biggl(\,\prod_{1\leq i<j\leq n}
\frac{u_i-tu_j}{u_i-u_j} \biggr)
=\frac{(t;t)_n}{(1-t)^n},
\end{equation}
we obtain
\[
\frac{(t;t)_{n-m}}{(1-t)^{n-m}}
\sum_{\substack{l_1,\dots,l_m=1\\ l_i\neq l_j}}^n
\prod_{i=1}^m\frac{y_{l_i}}{x_i-ty_{l_i}} 
\prod_{i=1}^{n-m} 
\prod_{j=1}^m \frac{Y_i-ty_{l_j}}{Y_i-y_{l_j}} 
\prod_{1\leq i<j\leq m}
\frac{x_i-y_{l_j}}{x_i-ty_{l_j}}\cdot
\frac{y_{l_i}-ty_{l_i}}{y_{l_i}-y_{l_j}} .
\]

If we denote the expression on the right of \eqref{omegaSS}
by $\bar{\omega}(x,y;t)$, and use that
\[
\prod_{i=1}^{n-m} \prod_{j=1}^m \frac{Y_i-ty_{l_j}}{Y_i-y_{l_j}} =
\prod_{\substack{i=1\\ i\neq l_1,\dots,l_m}}^n \prod_{j=1}^m 
\frac{y_i-ty_{l_j}}{y_i-y_{l_j}},
\]
the above calculations imply that
\begin{multline*}
\bar{\omega}(x,y;t)=\kappa(t)
\sum_{\substack{l_1,\dots,l_m=1\\ l_i\neq l_j}}^n
\sum_{w\in\Symm_x} w\biggl(\;
\prod_{i=1}^m\frac{y_{l_i}}{x_i-ty_{l_i}} 
\prod_{\substack{i=1 \\ i\neq l_1,\dots,l_m}}^n 
\prod_{j=1}^m \frac{y_i-ty_{l_j}}{y_i-y_{l_j}} \\
\times
\prod_{1\leq i<j\leq m}
\frac{x_i-tx_j}{x_i-x_j}\cdot
\frac{x_i-y_{l_j}}{x_i-ty_{l_j}}\cdot
\frac{y_{l_i}-ty_{l_i}}{y_{l_i}-y_{l_j}}
\biggr),
\end{multline*}
where 
\[
\kappa(t)=\frac{t^{m(m-n)} (1-t)^{2m}}{(t;t)_m}.
\]

The expression for $\omega(x,y;t)$ given by \eqref{itt}
is also a sum over the $l_i$ but unfortunately the two summands 
do not equate and some further manipulations of the sums are 
required.

To proceed we apply
\begin{equation}\label{msum}
\sum_{\substack{l_1,\dots,l_m=1\\ l_i\neq l_j}}^n f(y_l)
=\sum_{1\leq l_1<\cdots<l_m\leq n}\:
\sum_{w\in\Symm_{y_l}}w\bigl(f(y_l)\bigr),
\end{equation}
with $y_l=(y_{l_1},\dots,y_{l_m})$. Therefore
\begin{multline*}
\bar{\omega}(x,y;t)=
\kappa(t)
\sum_{1\leq l_1<\cdots<l_m\leq n}\: 
\sum_{w\in\Symm_x\times\Symm_{y_l}}
w \biggl(\;
\prod_{i=1}^m\frac{y_{l_i}}{x_i-ty_{l_i}} 
\prod_{\substack{i=1 \\ i\neq l_1,\dots,l_m}}^n 
\prod_{j=1}^m \frac{y_i-ty_{l_j}}{y_i-y_{l_j}} \\
\times
\prod_{1\leq i<j\leq m}
\frac{x_i-tx_j}{x_i-x_j}\cdot
\frac{x_i-y_{l_j}}{x_i-ty_{l_j}}\cdot
\frac{y_{l_i}-ty_{l_i}}{y_{l_i}-y_{l_j}}\biggr).
\end{multline*}

We now invoke the following lemma, which reduces to
\eqref{M14} for $v=u$.
\begin{lemma}\label{lemS2}
For $u=(u_1,\dots,u_n)$ and $v=(v_1,\dots,v_n)$ there holds
\begin{multline*}
\sum_{w\in\Symm_u\times\Symm_v} w\biggl(\;
\prod_{i=1}^n\frac{1}{u_i-tv_i} \prod_{1\leq i<j\leq n}
\frac{u_i-tu_j}{u_i-u_j}\cdot
\frac{u_i-v_j}{u_i-tv_j}\cdot
\frac{v_i-tv_j}{v_i-v_j} \biggr) \\
=\frac{(t;t)_n}{(1-t)^n}\sum_{w\in\Symm_v} w\biggl(\;
\prod_{i=1}^n \frac{1}{u_i-tv_i} \prod_{1\leq i<j\leq n}
\frac{u_i-v_j}{u_i-tv_j}\cdot
\frac{v_i-tv_j}{v_i-v_j}\biggr).
\end{multline*}
\end{lemma}
Since
\[
\prod_{i=1}^m y_{l_i}
\prod_{\substack{i=1 \\ i\neq l_1,\dots,l_m}}^n 
\prod_{j=1}^m \frac{y_i-ty_{l_j}}{y_i-y_{l_j}} 
\]
is symmetric in $y_l$, Lemma~\ref{lemS2} (with $(n,u,v)\to (m,x,y_l)$)
may be applied to yield
\begin{multline*}
\bar{\omega}(x,y;t)=
\kappa(t)\,\frac{(t;t)_m}{(1-t)^m}
\sum_{1\leq l_1<\cdots<l_m\leq n} \sum_{w\in\Symm_{y_l}}
w \biggl(\;
\prod_{i=1}^m \frac{y_{l_i}}{x_i-ty_{l_i}} \\
\times
\prod_{\substack{i=1 \\ i\neq l_1,\dots,l_m}}^n 
\prod_{j=1}^m \frac{y_i-ty_{l_j}}{y_i-y_{l_j}}
\prod_{1\leq i<j\leq m}
\frac{x_i-y_{l_j}}{x_i-ty_{l_j}} \cdot 
\frac{y_{l_i}-ty_{l_j}}{y_{l_i}-y_{l_j}} \biggr).
\end{multline*}
Reversing \eqref{msum} we finally get
\[
\bar{\omega}(x,y;t)=
\kappa(t)\,\frac{(t;t)_m}{(1-t)^m}
\sum_{\substack{l_1,\dots,l_m=1 \\ l_i\neq l_j}}^n
\prod_{i=1}^m \frac{y_{l_i}}{x_i-ty_{l_i}} 
\prod_{\substack{i=1 \\ i\neq l_1,\dots,l_m}}^n 
\prod_{j=1}^m \frac{y_i-ty_{l_j}}{y_i-y_{l_j}} 
\prod_{1\leq i<j\leq m}
\frac{x_i-y_{l_j}}{x_i-ty_{l_j}} \cdot 
\frac{y_{l_i}-ty_{l_j}}{y_{l_i}-y_{l_j}}.
\]
Comparing this with \eqref{itt} we see that $\bar{\omega}=\omega$
and the proof is complete except for a proof of Lemma~\ref{lemS2}.
\end{proof}

\begin{proof}[Proof of Lemma~\ref{lemS2}]
Defining
\begin{equation*}
g(u,v;t)=\prod_{i=1}^n \frac{1}{u_i-tv_i} \prod_{1\leq i<j\leq n}
\frac{u_i-v_j}{u_i-tv_j}\cdot\frac{v_i-tv_j}{v_i-v_j},
\end{equation*}
the proposition states that
\begin{equation}\label{xysym}
\sum_{w\in \Symm_u\times\Symm_v} w\biggl( g(u,v;t)
\prod_{1\leq i<j\leq n} \frac{u_i-tu_j}{u_i-u_j} \biggr)
=\frac{(t;t)_n}{(1-t)^n}
\sum_{w\in\Symm_v}w\bigl(g(u,v;t)\bigr).
\end{equation}
The difficulty is that it is unclear that the right-hand side
is symmetric in $u$. For example, when $n=2$ it reads
(without the $(u,v)$-independent prefactor)
\[
\frac{1}{u_1-tv_1}\cdot
\frac{1}{u_2-tv_2}\cdot
\frac{u_1-v_2}{u_1-tv_2}\cdot
\frac{v_1-tv_2}{v_1-tv_2} 
+\frac{1}{u_1-tv_2}\cdot
\frac{1}{u_2-tv_1}\cdot
\frac{u_1-v_1}{u_1-tv_1}\cdot
\frac{v_2-tv_1}{v_2-tv_1},
\]
which appears symmetric in $v$ only, but is in fact equal to
\begin{equation*}
\frac{(1+t) (tv_1 v_2+u_1 u_2) - t(v_1+v_2)(u_1+u_2)}
{(u_1-t v_1)(u_1-t v_2)(u_2-t v_1)(u_2-t v_2)}.
\end{equation*}

\medskip

Let $T_{k,u}\in \Symm_u$ by the $k$th adjacent transposition acting on $u$:
\begin{equation*}
T_{k,u}\bigl( f(u) \bigr)=
f(u_1,\dots,u_{k-1},u_{k+1},u_k,u_{k+2},\dots,u_n).
\end{equation*}
The $T_{k,u}$ for $1\leq k\leq n-1$ generate $\Symm_u$, and to prove
that the right-hand side of \eqref{xysym} is symmetric in $u$
it suffices to show that it is invariant under the action of the
$T_{k,u}$. That is, we must show that
\begin{equation*}
T_{k,u}\Bigl(\:\sum_{w\in \Symm_v}w\bigl(g(u,v;t)\bigr)\Bigr)=
\sum_{w\in\Symm_v}w\bigl(g(u,v;t)\bigr)
\end{equation*}
or, equivalently,
\begin{equation}\label{MS}
\sum_{w\in\Symm_v}w
\Bigl(\,T_{k,u}\bigl(g(u,v;t)\bigr)\Bigr)=
\sum_{w\in\Symm_v}w\bigl(g(u,v;t)\bigr)
\end{equation}
since $T_{k,u}$ commutes with the $v$-symmetrization.

A direct computation shows that
\begin{equation*}
T_{k,u} \bigl(g(u,v;t) \bigr)
=g(u,v;t)
-\frac{(u_k-u_{k+1})(v_{k+1}-t v_k)}{(u_k-v_{k+1})(u_{k+1}-t v_k)}\:
g(u,v;t).
\end{equation*}
Acting with $\Symm_v$ it thus follows that \eqref{MS} holds if
\begin{equation}\label{nul}
\sum_{w\in\Symm_v}w\bigl( h(u,v;t)\bigr)=0
\end{equation}
for
\begin{equation*}
h(u,v;t)=\frac{(v_{k+1}-t v_k)}{(u_k-v_{k+1})(u_{k+1}-t v_k)}\:
g(u,v;t).
\end{equation*}
Given an arbitrary permutation $w=(w_1,\dots,w_n)\in\Symm_v$ let
$w'\in\Symm_v$ be given by
\begin{equation*}
w'=(w_1,\dots,w_{k-1},w_{k+1},w_k,w_{k+2},\dots,w_n).
\end{equation*}
Another direct computation shows that 
\begin{equation*}
w\bigl(h(u,v;t)\bigr)=-w'\bigl(h(u,v;t)\bigr).
\end{equation*}
Therefore
\begin{equation*}
\sum_{w\in\Symm_v}w\bigl(h(u,v;t)\bigr)=
-\sum_{w\in\Symm_v}w\bigl(h(u,v;t)\bigr)
\end{equation*}
from which \eqref{nul} follows.

Now that the $u$-symmetry of the right-hand side of \eqref{xysym} 
has been established the rest is easy.
By \eqref{M14}
\begin{align*}
\text{RHS}\eqref{xysym}&=
\sum_{w\in\Symm_u}w\biggl(\, \frac{u_i-tu_j}{u_i-u_j} \biggr)
\sum_{w\in\Symm_v}w\bigl(g(u,v;t)\bigr)  \\
&=\sum_{w\in\Symm_u\times\Symm_v}w\biggl(\, \frac{u_i-tu_j}{u_i-u_j} \:
g(u,v;t) \biggl)=\text{LHS}\eqref{xysym}
\end{align*}
completing the proof.
\end{proof}

\begin{proof}[Proof of Lemma~\ref{lemdet}]
The entries of the determinant may be expanded by
\begin{equation*}
\frac{1}{(1-xy)(1-txy)}=
\sum_{\alpha=0}^{\infty}\qbin{\alpha+1}{\alpha}_t(xy)^{\alpha},
\end{equation*}
where
\[
\qbin{N}{k}_q=\frac{(q^{N-k+1};q)_k}{(q;q)_k}
\]
is a $q$-binomial coefficient.
By multilinearity this gives
\[
\det_{1\leq i,j\leq n}\biggl(\dots \biggr)
=\sum_{\alpha_1,\dots,\alpha_n=0}^{\infty}
\det_{1\leq i,j\leq n}\bigl( y_j^{\alpha_i} \bigr)\: x^{\alpha}
\qbin{\alpha+1}{\alpha}_t ,
\]
where 
\[\qbin{\alpha+1}{\alpha}_t =\prod_{i=1}^n 
\qbin{\alpha_i+1}{\alpha_i}_t.
\]
Since the summand vanishes when two (or more) of the summation indices
coincide and since the product of $t$-binomials is symmetric in 
$\alpha$, this may be rewritten as
\begin{align*}
\det_{1\leq i,j\leq n}\biggl(\dots \biggr)
&=\sum_{\alpha_1>\dots>\alpha_n=0}^{\infty}
\sum_{w\in \Symm_{\alpha}}
\det_{1\leq i,j\leq n}\bigl( y_j^{\alpha_{w_i}} \bigr)\: x^{w(\alpha)}
\qbin{\alpha+1}{\alpha}_t \\
&=\sum_{\alpha_1>\dots>\alpha_n=0}^{\infty}
\det_{1\leq i,j\leq n}\bigl( y_j^{\alpha_i} \bigr)\:
\qbin{\alpha+1}{\alpha}_t
\sum_{w\in \Symm_{\alpha}} \epsilon(w) x^{w(\alpha)} \\
&=\sum_{\alpha_1>\dots>\alpha_n=0}^{\infty}
\det_{1\leq i,j\leq n}\bigl( y_j^{\alpha_i} \bigr)
\det_{1\leq i,j\leq n}\bigl( x_j^{\alpha_i} \bigr)\:
\qbin{\alpha+1}{\alpha}_t
\end{align*}
where $\epsilon(w)$ in the second line denotes
the signature of the permutation $w$.

Setting $\alpha_i=\la_i+n-i+1$ and using \eqref{Schur}
this becomes
\[
\det_{1\leq i,j\leq n}\biggl(\dots \biggr)
=\Delta(x)\Delta(y)\sum_{\la} s_{\la}(x)s_{\la}(y)
\prod_{i=1}^n \qbin{\la_i+n-i+1}{\la_i+n-i}_t.
\]

Recalling that $m=n$ we have
\[\prod_{i=1}^n \qbin{\la_i+n-i+1}{\la_i+n-i}_t=\frac{g_{\la}(t;t,t)}{(1-t)^n}\]
so that
\[
\sum_{\la} g_{\la}(t;t,t) s_{\la}(x)s_{\la}(y)=
(1-t)^n \det_{1\leq i,j\leq n}\biggl(\dots \biggr)
\frac{1}{\Delta(x)\Delta(y)}.
\]
By \eqref{qist} the left-hand side may be recognised as the left-hand side of
\eqref{Fdef} for $m=n$, $q=t$ and $u=t$. Hence it may be replaced by the
corresponding right-hand side, leading to
\[
F(t;x,y;t)=
(1-t)^n \det_{1\leq i,j\leq n}\biggl(\dots \biggr)
\frac{\prod_{i,j=1}^n (1-x_iy_j)}{\Delta(x)\Delta(y)}
\]
as claimed by the lemma.
\end{proof}

\section{An identity for $q,t$-Littlewood--Richardson coefficients}
\label{secqtLR}
In our proof of the $\gsl_3$ $q$-binomial theorem \eqref{qbthmsl3}
we require the following identity for the 
$q,t$-Littlewood--Richardson coefficients.
\begin{theorem}\label{thm1}
Given integers $0\leq m\leq n$, let $\la$ and $\mu$ be 
partitions such that $l(\la)\leq m$ and $l(\mu)\leq n$.
Then
\begin{multline*}
\sum_{\omega,\nu} 
t^{n(\nu)-\abs{\omega}}f_{\omega\nu}^{\la}(q,t) 
V_{\nu,0}(u,1;q,t) u_0^{(n-m)}(P_{\mu/\omega})\,
\frac{(qt^{m-n-1};q,t)_{\nu}}{c'_{\nu}(q,t)} \\
=t^{n(\la)-m\abs{\mu}} V_{\la\mu}(u,1;q,t) 
u_0^{(n)}(P_{\mu})\,
\frac{(qt^{m-1};q,t)_{\la}}{c'_{\la}(q,t)}
\prod_{i=1}^m \prod_{j=1}^n 
\frac{(qt^{j-i+m-n-1};q)_{\la_i-\mu_j}}
{(qt^{j-i+m-n};q)_{\la_i-\mu_j}}.
\end{multline*}
\end{theorem}

Since $f_{\omega\nu}^{\la}(q,t)=0$ if $\omega\not\subseteq\la$ and
$P_{\mu/\omega}=0$ if $\omega\not\subseteq\mu$ we may add the
restrictions $\omega\subseteq\la$ and $\omega\subseteq\mu$ to
the sum over $\omega$.
It may in fact also be shown that the summand on the left vanishes
unless
\begin{equation}\label{van}
\la_i\geq \mu_{i+n-m} \quad\text{for~~$1\leq i\leq m$.}
\end{equation}
In other words, if $\mu^{\ast}$ is the partition formed by the
last $m$ parts of $\mu$ (i.e., $\mu^{\ast}=(\mu_{n-m+1},\dots,\mu_n)$)
then the summand vanishes unless $\mu^{\ast}\subseteq\la$.

To see this we recall from \cite[Equation (VI.7.13$'$)]{Macdonald95} that
\begin{equation*}
P_{\mu/\omega}(x_1,\dots,x_{n-m};q,t)=\sum_{T}\psi_T(q,t)x^T,
\end{equation*}
where the sum is over all semistandard Young tableaux $T$ of 
skew shape $\mu-\omega$ over the alphabet $\{1,\dots,n-m\}$;
$x^T$ is the monomial defined by $T$ and $\psi_T\in \F$.
For the shape $\mu-\omega$ to have an admissible filling it
must have at most $n-m$ boxes in each of its columns.
Hence $\omega_i\geq\mu_{i+n-m}$ for $1\leq i\leq m$.
Since we already established that the summand vanishes unless
$\omega\subseteq\la$, a necessary condition for nonvanishing of
the summand is thus given by \eqref{van}.
Since $1/(q;q)_{-N}=0$ for $N$ a positive integer,
it is easily seen that also the double product
on the right-hand side of the theorem vanishes unless
\eqref{van} holds.

\begin{proof}[Proof of Theorem~\ref{thm1}]
We start with \eqref{Fdef} with $\la$ replaced by $\eta$
and apply the homomorphisms 
$u_{\la;z}^{(m)}$ (acting on $x$) and $\displaystyle 
u_{\mu}^{(n)}$ (acting on $y$).
Using the homogeneity \eqref{hom} of the Macdonald polynomials
and recalling \eqref{Wlamudef} this leads to
\begin{multline}\label{eta}
\sum_{\eta} 
z^{\abs{\eta}} b_{\eta}(q,t)g_{\eta}(u;q,t)
u_{\la}^{(m)}(P_{\eta})u_{\mu}^{(n)}(P_{\eta}) \\
=W_{\la\mu}(u,z;q,t)\prod_{i=1}^m 
\frac{(zt^{n+m-i};q)_{\infty}}{(zt^{m-i};q)_{\infty}}
\prod_{i=1}^m \prod_{j=1}^n 
\frac{(zt^{n+m-i-j};q)_{\la_i+\mu_j}}{(zt^{n+m-i-j+1};q)_{\la_i+\mu_j}}.
\end{multline}

The summand on the left vanishes unless $l(\eta)\leq m$.
Assuming such $\eta$ we may twice use
the symmetry \eqref{symm} to rewrite the left-hand side as
\[
\text{LHS}\eqref{eta}=
\sum_{\eta} z^{\abs{\eta}} b_{\eta}(q,t) g_{\eta}(u;q,t)\,
\frac{u_{\eta}^{(m)}(P_{\la})u_{\eta}^{(n)}(P_{\mu}) 
u_0^{(m)}(P_{\eta})u_0^{(n)}(P_{\eta})}
{u_0^{(m)}(P_{\la})u_0^{(n)}(P_{\mu})}.
\]
Next we apply \eqref{skewdef2} as well as \eqref{hom} to get
\begin{align*}
u_{\eta}^{(n)}(P_{\mu})&=
P_{\mu}(q^{\eta_1}t^{n-1},\dots,q^{\eta_m}t^{n-m},t^{n-m-1},\dots,t,1;q,t) \\[2mm]
&=\sum_{\omega}P_{\omega}(q^{\eta_1}t^{n-1},\dots,q^{\eta_m}t^{n-m};q,t)
u_0^{(n-m)}(P_{\mu/\omega}) \\
&=\sum_{\omega}t^{(n-m)\abs{\omega}}u_{\eta}^{(m)}(P_{\omega})u_0^{(n-m)}(P_{\mu/\omega}).
\end{align*}
Thus
\[
\text{LHS}\eqref{eta}=\sum_{\eta,\omega} 
z^{\abs{\eta}} t^{(n-m)\abs{\omega}} b_{\eta}(q,t)
g_{\eta}(u;q,t) 
\:
\frac{u_0^{(n-m)}(P_{\mu/\omega})u_{\eta}^{(m)}(P_{\la}) u_{\eta}^{(m)}(P_{\omega})
u_0^{(m)}(P_{\eta})u_0^{(n)}(P_{\eta})}{u_0^{(m)}(P_{\la})u_0^{(n)}(P_{\mu})}.
\]
Next we use that
\begin{align*}
u_{\eta}^{(m)}(P_{\la}) u_{\eta}^{(m)}(P_{\omega})&=
u_{\eta}^{(m)}(P_{\la} \, P_{\omega}) \\[2mm]
&=u_{\eta}^{(m)}\Bigl(\, \sum_{\nu} f_{\omega\la}^{\nu}(q,t)P_{\nu}\Bigr) 
\qquad (\text{by \eqref{qtLR})} \\
&=\sum_{\nu} f_{\omega\la}^{\nu}(q,t)\, u_{\eta}^{(m)}(P_{\nu}) 
\end{align*}
to rewrite this as
\[
\text{LHS}\eqref{eta}=\sum_{\eta,\omega,\nu} 
z^{\abs{\eta}} t^{(n-m)\abs{\omega}}f_{\omega\la}^{\nu}(q,t) b_{\eta}(q,t) 
g_{\eta}(u;q,t) 
\: \frac{u_0^{(n-m)}(P_{\mu/\omega})u_{\eta}^{(m)}(P_{\nu}) 
u_0^{(m)}(P_{\eta})u_0^{(n)}(P_{\eta})}{u_0^{(m)}(P_{\la})u_0^{(n)}(P_{\mu})}.
\]
By one more application of \eqref{symm} this becomes
\[
\text{LHS}\eqref{eta}=\sum_{\eta,\omega,\nu} 
z^{\abs{\eta}} t^{(n-m)\abs{\omega}}f_{\omega\la}^{\nu}(q,t) b_{\eta}(q,t) 
g_{\eta}(u;q,t) 
\:
\frac{u_0^{(n-m)}(P_{\mu/\omega})
u_{\nu}^{(m)}(P_{\eta}) u_0^{(m)}(P_{\nu})
u_0^{(n)}(P_{\eta})}{u_0^{(m)}(P_{\la})u_0^{(n)}(P_{\mu})}.
\]

As a result of the previous manipulations the sum over $\eta$
corresponds to
\begin{align*}
\sum_{\eta} &
z^{\abs{\eta}} b_{\eta}(q,t) g_{\eta}(u;q,t)
u_{\nu}^{(m)}(P_{\eta}) u_0^{(n)}(P_{\eta}) &&\\
&=u_{\nu;z}^{(m)} u_0^{(n)}\biggl(
\sum_{\eta}
b_{\eta}(q,t) g_{\eta}(u;q,t)
P_{\eta}(x;q,t) P_{\eta}(y;q,t) \biggr) &&\\
&=u_{\nu;z}^{(m)} u_0^{(n)}\biggl(
F(u;x,y;t) \prod_{i=1}^m \prod_{j=1}^n 
\frac{(tx_i y_j;q)_{\infty}}{(x_i y_j;q)_{\infty}} \biggr) 
&&(\text{by \eqref{Fdef}}) \\
&=u_{\nu;z}^{(m)} u_0^{(n)}\Bigl(F(u;x,y;t)\Bigr)
\frac{(zt^{m-1};q,t)_{\nu}}{(zt^{n+m-1};q,t)_{\nu}}
\prod_{i=1}^m \frac{(zt^{n+m-i};q)_{\infty}} {(zt^{m-i};q)_{\infty}} 
&& \\
&=W_{\nu,0}(u,z;q,t)\,
\frac{(zt^{m-1};q,t)_{\nu}}{(zt^{n+m-1};q,t)_{\nu}}
\prod_{i=1}^m \frac{(zt^{n+m-i};q)_{\infty}} {(zt^{m-i};q)_{\infty}}
&&(\text{by \eqref{Wlamudef}).}
\end{align*}
We thus arrive at
\begin{multline*}
\text{LHS}\eqref{eta}=
\prod_{i=1}^m \frac{(zt^{n+m-i};q)_{\infty}}{(zt^{m-i};q)_{\infty}} 
\\ \times
\sum_{\omega,\nu} 
t^{(n-m)\abs{\omega}}f_{\omega\la}^{\nu}(q,t)\, W_{\nu,0}(u,z;q,t)\,
\frac{u_0^{(n-m)}(P_{\mu/\omega})u_0^{(m)}(P_{\nu})}
{u_0^{(n)}(P_{\mu})u_0^{(m)}(P_{\la})}\:
\frac{(zt^{m-1};q,t)_{\nu}}{(zt^{n+m-1};q,t)_{\nu}}.
\end{multline*}
Finally equating this with the right-hand side of \eqref{eta} yields
\begin{multline}\label{WW}
\sum_{\omega,\nu} 
t^{(n-m)\abs{\omega}}f_{\omega\la}^{\nu}(q,t)\, W_{\nu,0}(u,z;q,t)\:
\frac{u_0^{(n-m)}(P_{\mu/\omega}) u_0^{(m)}(P_{\nu})}
{u_0^{(n)}(P_{\mu})u_0^{(m)}(P_{\la})}\,
\frac{(zt^{m-1};q,t)_{\nu}}{(zt^{n+m-1};q,t)_{\nu}} \\
=W_{\la\mu}(u,z;q,t)
\prod_{i=1}^m \prod_{j=1}^n 
\frac{(zt^{n+m-i-j};q)_{\la_i+\mu_j}}{(zt^{n+m-i-j+1};q)_{\la_i+\mu_j}}.
\end{multline}
Both sides of this identity trivially vanish if
$l(\la)>m$. Furthermore, the summand on the left
vanishes if $l(\nu)>m$.
Hence we may without loss of generality assume in
the following that $l(\la)\leq m$ and 
$l(\nu)\leq m$. (The latter of course refers to
a restriction on the summation index.)  
We may also assume that the largest part of $\nu$ is bounded since
$f_{\omega\la}^{\nu}=0$ if $\abs{\omega}+\abs{\la}\neq\abs{\nu}$
and $P_{\mu/\omega}=0$ if $\omega\not\subseteq\mu$.
In particular $\nu_1\leq \abs{\la}+\abs{\mu}$.

The above considerations imply that
$\la,\nu\subseteq (N^m)$ for sufficiently large $N$.
Given such $N$ we can define the partitions 
$\hat{\la}$ and $\hat{\mu}$ as
the complements of $\la$ and $\nu$ with
respect to $(N^m)$, i.e.,
$\hat{\la}_i=N-\la_{m+1-i}$ and
$\hat{\nu}_i=N-\nu_{m+1-i}$ for $1\leq i\leq m$.

We now replace $z\to q^{1-m-N}/z$, $\la\to\hat{\la}$ and $\nu\to\hat{\nu}$
in \eqref{WW}, and then eliminate the hats.
For this we need the easily established
\[
W_{\hat{\la}\mu}(u,q^{1-m-N}/z;q,t)
=V_{\la\mu}(u,z;q,t)
\]
as well as \cite[page 263]{Warnaar}
\begin{equation*}
f_{\omega\hat{\la}}^{\hat{\nu}}(q,t)=
t^{n(\nu)-n(\la)}
f_{\omega\nu}^{\la}(q,t)\,
\frac{(qt^{m-1};q,t)_{\nu}}{(qt^{m-1};q,t)_{\la}}\,
\frac{c'_{\la}(q,t)}{c'_{\nu}(q,t)}\,
\frac{u_0^{(m)}(P_{\la})}{u_0^{(m)}(P_{\nu})},
\end{equation*}
\cite[Equation (4.1)]{BF99}
\begin{equation*}
\frac{(a;q,t)_{\hat{\la}}}{(b;q,t)_{\hat{\la}}}
=\Bigl(\frac{b}{a}\Bigr)^{\abs{\la}}
\frac{(a;q,t)_{(N^m)}}{(b;q,t)_{(N^m)}}\,
\frac{(q^{1-N}t^{m-1}/b;q,t)_{\la}}
{(q^{1-N}t^{m-1}/a;q,t)_{\la}},
\end{equation*}
and
\begin{equation*}
u_0^{(m)}(P_{\hat{\la}})=t^{\binom{m}{2}N+(1-m)\abs{\la}}
u_0^{(m)}(P_{\la}).
\end{equation*}
This last result follows from \cite[Equation (4.3)]{BF99}
\begin{equation*}
P_{\hat{\la}}(x;q,t)=(x_1\cdots x_m)^N\, P_{\la}(x^{-1};q,t)
\end{equation*}
and the homogeneity \eqref{hom}.
As a result we arrive at
\begin{multline*}
\sum_{\omega,\nu} 
t^{n(\nu)-\abs{\omega}}f_{\omega\nu}^{\la}(q,t) 
V_{\nu,0}(u,z;q,t) 
u_0^{(n-m)}(P_{\mu/\omega})\,
\frac{(qt^{m-1},zqt^{m-n-1};q,t)_{\nu}}
{c'_{\nu}(q,t)\, (zqt^{m-1};q,t)_{\nu}} \\
=t^{n(\la)-m\abs{\mu}} 
V_{\la\mu}(u,z;q,t) 
u_0^{(n)}(P_{\mu}) 
\:
\frac{(qt^{m-1};q,t)_{\la}}{c'_{\la}(q,t)}
\prod_{i=1}^m \prod_{j=1}^n 
\frac{(zqt^{j-i+m-n-1};q)_{\la_i-\mu_j}}
{(zqt^{j-i+m-n};q)_{\la_i-\mu_j}},
\end{multline*}
where we have also that
$f_{\omega\nu}^{\la}=0$ if $\abs{\omega}+\abs{\nu}\neq \abs{\la}$,
and
\begin{equation*}
\frac{(a;q)_{N-k}}{(b;q)_{N-k}}=
\frac{(a;q)_N}{(b;q)_N}\,
\frac{(q^{1-N}/b;q)_k}{(q^{1-N}/a;q)_k}\Bigl(\frac{b}{a}\Bigr)^k.
\end{equation*}
Finally specializing $z=1$ complete the proof.
\end{proof}

\section{$\gsl_3$ basic hypergeometric series}\label{SecA2}

Below we will give two different definitions of $\gsl_3$ basic 
hypergeometric series, denoted Type I and Type II respectively.
To cover both types at once we introduce the function
$V_{\la\mu}(q,t)$ which is either given by
\begin{align*}
V_{\la\mu}(q,t)&=V_{\la\mu}(1,1;q,t) &&\text{Type I}\\
\intertext{or by}
V_{\la\mu}(q,t)&=q^{-\abs{\la}}V_{\la\mu}(t^{n-m+1},1;q,t)
&&\text{Type II}.
\end{align*}
Note that it follows from \eqref{Vlamudef} and \eqref{omegaF}
that for Type I series, 
\[
V_{\la\mu}(q,t)=u_{\la}^{(m)}u_{\mu}^{(n)}
\bigl(\omega(x,y;t)\bigr).
\]

{}From \eqref{VPS} and \eqref{VPS2} we see that
regardless of our choice of $V_{\la\mu}(q,t)$
\begin{equation}\label{Vprod}
V_{\la,0}(q,t)=\prod_{i=1}^m \frac{1-t^{m-n-i}}{1-q^{\la_i}t^{m-n-i}}=
\frac{(t^{m-n-1};q,t)_{\la}}{(qt^{m-n-1};q,t)_{\la}}.
\end{equation}

It is important to observe that $V_{\la\mu}(q,t)$ 
does not merely depend of the partitions $\la$ and $\mu$
but also on the integers $m$ and $n$. (We tacitly assume
that $l(\la)\leq m$ and $l(\mu)\leq n$.)
These integers are mostly assumed to be fixed, 
but occasionally we will relate series labelled
by $(m,n)$ to those labelled by $(m-1,n-1)$.
If we write
$V^{(m,n)}_{\la\mu}(q,t)$ instead of $V_{\la\mu}(q,t)$
it follows from Lemma~\ref{lemS1} that $V^{(m,n)}_{\la\mu}(q,t)$
only depends on the difference $n-m$. Specifically,
\begin{equation}\label{Vmn}
V^{(m,n)}_{\la\mu}(q,t)=V^{(m-1,n-1)}_{\la\mu}(q,t)
\end{equation}
provided of course that $l(\la)\leq m-1$ and $l(\mu)\leq n-1$.

To reduce the length of many of the subsequent formulae we introduce
another rational function $\Omega_{\la\mu}(q,t)$ as
\begin{equation}\label{Omegadef}
\Omega_{\la\mu}(q,t)=
V_{\la\mu}(q,t) \, (qt^{m-1};q,t)_{\la}
\prod_{i=1}^m\prod_{j=1}^n
\frac{(q t^{j-i+m-n-1};q)_{\la_i-\mu_j}}
{(q t^{j-i+m-n};q)_{\la_i-\mu_j}},
\end{equation}
where $\la$ and $\mu$ are partitions such that
$l(\la)\leq m$ and $l(\mu)\leq n$.

Two easily established results for $\Omega_{\la\mu}(q,t)$ are
\begin{equation}\label{Omega0}
\Omega_{\la,0}(q,t)=(t^{m-n-1};q)_{\la}
\end{equation}
and, displaying the $(m,n)$ dependence,
\begin{equation}\label{Omegamn}
\Omega_{\la\mu}^{(m,n)}(q,t)=
\Omega_{\la\mu}^{(m-1,n-1)}(q,t)\, t^{\abs{\mu}}
\frac{(t^{n-1};q,t)_{\mu}}{(t^n;q,t)_{\mu}}
\end{equation}
for $l(\la)\leq m-1$ and $l(\mu)\leq n-1$.
Equation \eqref{Omega0} follows from \eqref{Vprod} and
\[
\prod_{i=1}^m\prod_{j=1}^n
\frac{(q t^{j-i+m-n-1};q)_{\la_i-\mu_j}}
{(q t^{j-i+m-n};q)_{\la_i-\mu_j}}\bigg|_{\mu=0}
=\frac{(qt^{m-n-1};q,t)_{\la}}{(qt^{m-1};q,t)_{\la}},
\]
and \eqref{Omegamn} follows from \eqref{Vmn} and
\[
\prod_{i=1}^m\prod_{j=1}^n
\frac{(q t^{j-i+m-n-1};q)_{\la_i-\mu_j}}
{(q t^{j-i+m-n};q)_{\la_i-\mu_j}}\bigg|_{\la_m=\mu_n=0}
=t^{\abs{\mu}}\,
\frac{(qt^{m-2};q,t)_{\la}}{(qt^{m-1};q,t)_{\la}}\:
\frac{(t^{n-1};q,t)_{\mu}}{t^n;q,t)_{\mu}}.
\]

We can now state the main definition of this section.
\begin{definition}[$\gsl_3$ basic hypergeometric series]\label{defl}
Let $x=(x_1,\dots,x_m)$ and $y=(y_1,\dots,y_n)$ such that
$0\leq m\leq n$. Then
\begin{multline}\label{seriesdef}
{_{r+1}\Phi_r}\hyp{a_1,\dots,a_{r+1}}{b_1,\dots,b_r}{q,t;x,y}
=\prod_{i=1}^m 
\frac{(x_i;q)_{\infty}}{(x_it^{m-n-1};q)_{\infty}} \\
\times\sum_{\la,\mu} t^{n(\la)+n(\mu)}\,
\frac{P_{\la}(x;q,t)}{c'_{\la}(q,t)}\,
\frac{P_{\mu}(y;q,t)}{c'_{\mu}(q,t)}\,
\frac{(a_1,\dots,a_{r+1};q,t)_{\mu}}{(b_1,\dots, b_r;q,t)_{\mu}} 
\, \Omega_{\la\mu}(q,t), 
\end{multline}
where the sum is over partitions $\la$ and $\mu$ such that
$l(\la)\leq m$, $l(\mu)\leq n$ and 
\begin{equation}\label{order}
\la_i\geq \mu_{i-m+n} \quad \text{for~~$1\leq i\leq m$.}
\end{equation}
\end{definition}

\noindent\textbf{Remarks.}
\begin{enumerate}
\item
The restrictions on the sum may alternatively be expressed by 
the inequalities \cite[Equation (2.4)]{TV03}
\begin{equation*}
\begin{array}{cccccccccccccc}
&&&&\la_1&\GEQ &\la_2 &\GEQ&\dots&\GEQ&\la_m&&& \\[-1mm]
&&&&\text{\begin{turn}{270}$\geq$\end{turn}}&&
\text{\begin{turn}{270}$\geq$\end{turn}} &&&&
\text{\begin{turn}{270}$\geq$\end{turn}} &&& \\[3mm]
\mu_1&\GEQ&\cdots&\GEQ &\mu_{n-m+1}&\GEQ &\mu_{n-m+2}&\GEQ&\dots&\GEQ
&\mu_n&\GEQ & 0
\end{array}
\end{equation*}
\item
The prefactor 
\begin{equation*}
\prod_{i=1}^m \frac{(x_i;q)_{\infty}}{(x_it^{m-n-1};q)_{\infty}}
\end{equation*}
in the definition has been included to simplify subsequent
formulae, and implies that for $a_1=1$ 
the $\gsl_3$ series simplifies to $1$, see Lemma~\ref{lemma1} below.
\item
The main reason for attaching the label $\gsl_3$ to the series of
Definition~\ref{defl} is the connection with the $\gsl_3$
discrete exponential and continous Selberg integrals of 
Tarasov and Varchenko, see page~\pageref{page} for details.

We should also mention that we are
currently developing a theory of $\gsl_n$ 
basic hypergeometric series \cite{Warnaarb}. 
In such series, a Macdonald polynomial is attached to each
vertex of the $\gsl_n$ Dynkin diagram, and the corresponding
$\gsl_n$ $q$-binomial theorem may be expressed concisely 
in terms of the data of the undelying Lie algebra.
\item
Finally we remark that nearly all our results 
involve non-terminating $\gsl_3$ series. To ensure
convergence we implicitly assume that
\begin{equation*}
\max\{\abs{q},\abs{t},\abs{x_1},\dots,
\abs{x_m},\abs{y_1},\dots,\abs{y_n}\}<1
\end{equation*}
whenever necessary.
\end{enumerate}

Our most important results for $\gsl_3$ basic
hypergeometric series are two generalizations of the $q$-binomial 
theorem. First however, we state several elementary properties of 
the series.
In all of the results below the parameters 
$a_1,\dots,a_{r+1}$ and $b_1,\dots,b_r$
act as dummies, and to shorten some of the equations
we abbreviate these sequences by $A$ and $B$ respectively.
\begin{lemma}\label{lem0}
We have
\begin{equation*}
{_{r+1}\Phi_r}\hyp{A}{B}{q,t;x,(0^n)}=1.
\end{equation*}
\end{lemma}

\begin{proof}[Proof of Lemma~\ref{lem0}]
Since $P_{\mu}((0^n);q,t)=\delta_{\mu,0}$
we get
\[
{_{r+1}\Phi_r}\hyp{A}{B}{q,t;x,(0^n)}
=\prod_{i=1}^m \frac{(x_i;q)_{\infty}}{(x_it^{m-n-1};q)_{\infty}}
\sum_{\la} t^{n(\la)}
\frac{P_{\la}(x;q,t)}{c'_{\la}(q,t)}\, \Omega_{\la,0}(q,t).
\]
Thanks to \eqref{Omega0} this is
\[
{_{r+1}\Phi_r}\hyp{A}{B}{q,t;x,(0^n)}
={_1\Phi_0}\hyp{t^{m-n-1}}{\textbf{--}}{q,t;x}
\prod_{i=1}^m \frac{(x_i;q)_{\infty}}{(x_it^{m-n-1};q)_{\infty}},
\]
where on the right we have used definition \eqref{Phirs}
of the $\gsl_2$ Kaneko--Macdonald series.
Summing the $_1\Phi_0$ series by the $q$-binomial theorem 
\eqref{Phi10sum} results in the claim of the lemma.
\end{proof}

\begin{lemma}\label{lemma1}
We have
\begin{equation*}
{_{r+1}\Phi_r}\hyp{1,a_2,\dots,a_{r+1}}{b_1,\dots,b_r}{q,t;x,y}=1.
\end{equation*}
\end{lemma}

\begin{proof}[Proof of Lemma~\ref{lemma1}]
When $a_1=1$ the summand vanishes unless $\mu=0$.
The proof is thus a repeat of the proof of Lemma~\ref{lem0}.
\end{proof}

The next two lemmas relate $\gsl_3$ 
series with labels $(n,m)$ and $(n-1,m-1)$. 
Recall the notation introduced in Section~\ref{secbiF}.
\begin{lemma}[(Stability 1)]\label{lemy0}
With $u_{0;z}^{(n)}$ acting on $y$ and
$u_{0;tz}^{(n-1)}$ acting on $y^{(n)}$, we have
\[
u_{0;z}^{(n)}\biggl({_{r+1}\Phi_r}\hyp{A}{B}{q,t;(x^{(m)},0),y}\biggr) 
=u_{0;tz}^{(n-1)}\biggl(
{_{r+1}\Phi_r}\hyp{A}{B}{q,t;x^{(m)},y^{(n)}}\biggr).
\]
\end{lemma}

\begin{lemma}[(Stability 2)]\label{lemxy0}
We have
\[
{_{r+1}\Phi_r}\hyp{A}{B}{q,t,(x^{(m)},0),(y^{(n)},0)}
={_{r+2}\Phi_{r+1}}\hyp{t^{n-1},A}{t^n,B}{q,t;x^{(m)},ty^{(n)}}.
\]
\end{lemma}

Iterating the two types of stability leads to
\begin{equation*}
u_{0;z}^{(n)}\biggl(
{_{r+1}\Phi_r}\hyp{A}{B}{q,t;(0^m),y}\biggr)
=u_{0;t^mz}^{(n-m)}\biggl(
{_{r+1}\Phi_r}\hyp{A}{B}{q,t;y^{(n-m+1,\dots,n)}}\biggr)
\end{equation*}
and
\[
{_{r+1}\Phi_r}\hyp{A}{B}{q,t;(0^m),(y^{(n-m+1,\dots,n)},0^m)}  
={_{r+2}\Phi_{r+1}}\hyp{t^{n-m},A}{t^n,B}{q,t;t^m y^{(n-m+1,\dots,n)}}.
\]
Note that both right-hand sides involve the
$\gsl_2$ Kaneko--Macdonald series.

\begin{proof}[Proof of Lemmas~\ref{lemy0} and \ref{lemxy0}]
Because we are comparing series for different $(m,n)$ values
we write $\Omega_{\la\mu}^{(m,n)}$ instead of $\Omega_{\la\mu}$.

If $x_m=0$ only partitions of length strictly less than $m$ 
contribute to the sum over $\la$. But if $\la_m=0$ then the 
inequality $0\leq \mu_n\leq\la_m$ implies that also $\mu_n=0$. 
Hence we may use \eqref{Omegamn} and the homogeneity of the 
Macdonald polynomials to obtain
\begin{multline*}
{_{r+1}\Phi_r}\hyp{A}{B}{q,t;(x^{(m)},0),y} 
=\prod_{i=1}^{m-1}\frac{(x_i;q)_{\infty}}{(x_it^{m-n-1};q)_{\infty}}
 \\ \times
\sum_{\la,\mu} t^{n(\la)+n(\mu)}
\frac{P_{\la}(x^{(m)};q,t)}{c'_{\la}(q,t)}\,
\frac{P_{\mu}(ty;q,t)}{c'_{\mu}(q,t)}\, 
\frac{(t^{n-1},A;q,t)_{\mu}}{(t^n,B;q,t)_{\mu}} \,
\Omega_{\la\mu}^{(m-1,n-1)}(q,t),
\end{multline*}
where the sum is over partitions $\la$ and $\mu$ such that
$l(\la)\leq m-1$, $l(\mu)\leq n-1$ and 
\begin{equation*}
\la_i\geq \mu_{i-m+n} \quad \text{for~~$1\leq i\leq m-1$.}
\end{equation*}
All terms on the right-hand side 
depend on $n-1$ and $m-1$ except for $P_{\mu}(ty;q,t)$,
since $y=(y_1,\dots,y_n)$. We can either make the obvious
choice $y_n=0$ and use the stability of the Macdonald polynomial:
$P_{\mu}(t(y^{(n)},0);q,t)=
P_{\mu}(ty^{(n)};q,t)$ to obtain
Lemma~\ref{lemxy0}, or we can specialize $y$.
In the latter case we may use that for $l(\mu)\leq n-1$
\begin{equation}\label{nnmin1}
u_{0;z}^{(n)}\bigl(P_{\mu}(y;q,t)\bigr)=
u_{0;z}^{(n-1)}\bigl(P_{\mu}(y^{(n)};q,t)\bigr)
\frac{(t^n;q,t)_{\mu}}{(t^{n-1};q,t)_{\mu}}
\end{equation}
as follows from \eqref{PS}. Therefore
\begin{multline*}
u_{0;z}^{(n)}\biggl({_{r+1}\Phi_r}\hyp{A}{B}{q,t;(x^{(m)},0),y}\biggr) 
=\prod_{i=1}^{m-1}
\frac{(x_i;q)_{\infty}}{(x_it^{m-n-1};q)_{\infty}}\\ \times
\sum_{\la,\mu} t^{n(\la)+n(\mu)}
\frac{P_{\la}(x^{(m)};q,t)}{c'_{\la}(q,t)}\, 
\frac{u_{0;tz}\bigl(P_{\mu}(y^{(n-1)};q,t)\bigr)}{c'_{\mu}(q,t)} \,
\frac{(A;q,t)_{\mu}}{(B;q,t)_{\mu}} \,
\Omega^{(m-1,n-1)}_{\la\mu}(q,t),
\end{multline*}
in accordance with the right-hand side of Lemma~\ref{lemy0}.
\end{proof}

Our next result implies all previous four lemmas,
but unlike the latter it is not elementary,
requiring Theorem~\ref{thm1} for its proof.
\begin{proposition}\label{prop}
Fix $\sigma$ as
\begin{equation}\label{sigma}
\sigma=\begin{cases}
0 &\text{for Type I} \\
1 &\text{for Type II}
\end{cases}
\end{equation}
and let $X=(X_1,\dots,X_n)$ be given by
\[
X_i=\begin{cases} 
q^{-\sigma}t^{-1}x_i &\text{for $1\leq i\leq m$} \\
t^{n-i}&\text{for $m+1\leq i\leq n$.}
\end{cases}
\]
Then
\begin{equation}\label{yX}
{_{r+1}\Phi_r}\hyp{A}{B}{q,t;x,y}
=\sum_{\mu}
t^{n(\mu)+m\abs{\mu}}\,
\frac{P_{\mu}(y;q,t)P_{\mu}(X;q,t)}{c'_{\mu}(q,t)\,u_0^{(n)}(P_{\mu})}\,
\frac{(A;q,t)_{\mu}}{(B;q,t)_{\mu}}.
\end{equation}
\end{proposition}
Note that by taking $y=(0^n)$ or $a_1=1$ the summand vanishes unless
$\mu=0$ leading to Lemmas~\ref{lem0} and \ref{lemma1}.
Also the Lemmas~\ref{lemy0} and \ref{lemxy0} immediately follow
from the proposition be it that the latter also requires
\eqref{nnmin1}.
For example, applying $u_{0;z}^{(n)}$ acting on $y$ to 
\eqref{yX} yields
\begin{equation*}
u_{0;z}^{(n)}\biggl(
{_{r+1}\Phi_r}\hyp{A}{B}{q,t;x,y}\biggr)
=\sum_{\mu}
z^{\abs{\mu}} t^{n(\mu)+m\abs{\mu}}\,
\frac{P_{\mu}(X;q,t)}{c'_{\mu}(q,t)}\,
\frac{(A;q,t)_{\mu}}{(B;q,t)_{\mu}}.
\end{equation*}
Not only does this make Lemma~\ref{lemy0} obvious but it in fact implies
the following more general (and more important) result.
\begin{corollary}\label{C1}
With the same notation as Proposition~\ref{prop} we have
\begin{equation}\label{C1eq}
u_{0;z}^{(n)}\biggl(
{_{r+1}\Phi_r}\hyp{A}{B}{q,t;x,y}\biggr)
={_{r+1}\Phi_r}\hyp{A}{B}{q,t;zt^mX}.
\end{equation}
\end{corollary}
Note that on the right we have the $\gsl_2$ Kaneko--Macdonald series.

There is another important corollary of Proposition~\ref{prop}.
If we take $m=n$ then 
\[
P_{\mu}(X;q,t)=q^{-\sigma\abs{\mu}}t^{-\abs{\mu}} P_{\mu}(x;q,t).
\]
Hence for $m=n$ the series \eqref{yX} is invariant under the interchange of $x$ and $y$.
\begin{corollary}\label{C2}
For $m=n$, i.e., $x=(x_1,\dots,x_n)$ and $y=(y_1,\dots,y_n)$, there holds
\begin{equation*}
{_{r+1}\Phi_r}\hyp{A}{B}{q,t;x,y}
={_{r+1}\Phi_r}\hyp{A}{B}{q,t;y,x}.
\end{equation*}
\end{corollary}

Using the above two corollaries it is straightforward to
prove several $q$-binomial theorems for $\gsl_3$ series.
First however we shall prove Proposition~\ref{prop}.
\begin{proof}[Proof of Proposition~\ref{prop}]
Recalling definition \eqref{Omegadef} and using \eqref{Vprod},
Theorem~\ref{thm1} may be rewritten as
\begin{multline*}
\Omega_{\la\mu}(q,t)=
\sum_{\omega,\nu} 
t^{n(\nu)-n(\la)+m\abs{\mu}-\abs{\omega}}f_{\omega\nu}^{\la}(q,t)\,
c'_{\la}(q,t)\\
\times \frac{V_{\nu,0}(u,1;q,t)}{V_{\nu,0}(q,t)} \,
\frac{V_{\la\mu}(q,t)}{V_{\la\mu}(u,1;q,t)}\,
\frac{u_0^{(n-m)}(P_{\mu/\omega})}{u_0^{(n)}(P_{\mu})}\,
\frac{(t^{m-n-1};q,t)_{\nu}}
{c'_{\nu}(q,t)}.
\end{multline*}
Taking $u=1$ or $u=t^{n-m+1}$, so that
\[
\frac{V_{\nu,0}(u,1;q,t)}{V_{\nu,0}(q,t)} \,
\frac{V_{\la\mu}(q,t)}{V_{\la\mu}(u,1;q,t)} \to
q^{-\sigma(\abs{\la}-\abs{\nu})},
\]
and using that $f_{\omega\nu}^{\la}=0$ if
$\abs{\omega}+\abs{\nu}\neq\abs{\la}$
we obtain
\[
\Omega_{\la\mu}(q,t)=\sum_{\omega,\nu} 
t^{n(\nu)-n(\la)+m\abs{\mu}-\abs{\omega}}q^{-\sigma\abs{\omega}}
f_{\omega\nu}^{\la}(q,t)\, c'_{\la}(q,t)
\: \frac{u_0^{(n-m)}(P_{\mu/\omega})}{u_0^{(n)}(P_{\mu})}\,
\frac{(t^{m-n-1};q,t)_{\nu}}{c'_{\nu}(q,t)}.
\]
Substituting this in the definition \eqref{seriesdef}
of the $\gsl_3$ basic hypergeometric series leads to
\begin{multline*}
{_{r+1}\Phi_r}\hyp{A}{B}{q,t;x,y}
=\prod_{i=1}^m
\frac{(x_i;q)_{\infty}}{(x_it^{m-n-1};q)_{\infty}} \\
\times\sum_{\la,\mu,\nu,\omega} 
t^{n(\mu)+n(\nu)+m\abs{\mu}-\abs{\omega}}q^{-\sigma \abs{\omega}}\,
\frac{P_{\mu}(y;q,t)}{c'_{\mu}(q,t)}\,
\frac{(A;q,t)_{\mu}}{(B;q,t)_{\mu}} \,
\frac{u_0^{(n-m)}(P_{\mu/\omega})}{u_0^{(n)}(P_{\mu})} \\
\times 
\frac{(t^{m-n-1};q,t)_{\nu}}{c'_{\nu}(q,t)}\,
f_{\omega\nu}^{\la}(q,t) 
P_{\la}(x;q,t).
\end{multline*}
Now performing the sum over $\la$ by \eqref{qtLR} yields
\begin{multline*}
{_{r+1}\Phi_r}\hyp{A}{B}{q,t;x,y}
=\prod_{i=1}^m
\frac{(x_i;q)_{\infty}}{(x_it^{m-n-1};q)_{\infty}} \\
\times\sum_{\mu,\nu,\omega} 
t^{n(\mu)+n(\nu)+m\abs{\mu}-\abs{\omega}}q^{-\sigma\abs{\omega}}\,
\frac{P_{\mu}(y;q,t)}{c'_{\mu}(q,t)}\,
\frac{(A;q,t)_{\mu}}{(B;q,t)_{\mu}} \,
\frac{u_0^{(n-m)}(P_{\mu/\omega})}{u_0^{(n)}(P_{\mu})} \\
\times 
\frac{(t^{m-n-1};q,t)_{\nu}}{c'_{\nu}(q,t)}\,
P_{\nu}(x;q,t) P_{\omega}(x;q,t).
\end{multline*}
The next simplification arises by noting that the sum over $\nu$
corresponds to a summable $\gsl_2$ Kaneko--Macdonald series:
\begin{equation*}
{_1\Phi_0}\hyp{qt^{m-n-1}}{\text{--}}{q,t;x}
=\prod_{i=1}^m \frac{(x_it^{m-n-1};q)_{\infty}}{(x_i;q)_{\infty}}
\end{equation*}
by \eqref{Phirs} and \eqref{Phi10sum}. Hence
\[
{_{r+1}\Phi_r}\hyp{A}{B}{q,t;x,y}
=\sum_{\mu,\omega}
t^{n(\mu)+m\abs{\mu}-\abs{\omega}}q^{-\sigma\abs{\omega}}\,
\frac{P_{\mu}(y;q,t)}{c'_{\mu}(q,t)}\,
\frac{(A;q,t)_{\mu}}{(B;q,t)_{\mu}}\, 
\frac{u_0^{(n-m)}(P_{\mu/\omega})}{u_0^{(n)}(P_{\mu})} \,
P_{\omega}(x;q,t).
\]
Next we use the homogeneity \eqref{hom} of $P_{\omega}$,
the definition \eqref{eval} of the principal specialization $u_0^{(n-m)}$
and the definition \eqref{skewdef2} of the skew Macdonald polynomials
to perform the sum over $\omega$;
\begin{align*}
\sum_{\omega}
(q^{\sigma}t)^{-\abs{\omega}} u_0^{(n-m)}(P_{\mu/\omega}) P_{\omega}(x;q,t)
&=\sum_{\omega}
u_0^{(n-m)}(P_{\mu/\omega}) P_{\omega}(q^{-\sigma}t^{-1}x;q,t) \\
&=P_{\mu}(X;q,t),
\end{align*}
where $X=(q^{-\sigma}t^{-1}x,t^{n-m-1},\dots,t,1)$.
The resulting identity is \eqref{yX}.
\end{proof}

{}From Corollary~\ref{C1} it is clear that whenever an
$\gsl_2$ series is summable this implies a
corresponding sum for $\gsl_3$ series.
The most obvious choice is to set $r=0$ in Corollary~\ref{C1}
so that the right-hand side of \eqref{C1eq}
may be summed by the Kaneko--Macdonald 
$q$-binomial theorem \eqref{Phi10sum}.
Hence
\begin{align*}
u_{0;z}^{(n)}\biggl({_1\Phi_0}\hyp{a}{\text{--}}{q,t;x,y}\biggr)
&=\prod_{i=1}^n \frac{(azt^mX_i;q)_{\infty}}{(zt^mX_i;q)_{\infty}} \\
&=\prod_{i=1}^m \frac{(azq^{-\sigma}t^{m-1}x_i;q)_{\infty}}
{(zq^{-\sigma}t^{m-1}x_i;q)_{\infty}}
\prod_{i=m+1}^n \frac{(azt^{m+n-i};q)_{\infty}}{(zt^{m+n-i};q)_{\infty}}.
\end{align*}

\begin{theorem}[(First $\gsl_3$ $q$-binomial theorem)]\label{FV}
For $x=(x_1,\dots,x_m)$ and $y=z(1,t,\dots,t^{n-1})$ we have
\[
{_1\Phi_0}\hyp{a}{\text{--}}{q,t;x,y}
=\prod_{i=1}^m \frac{(azt^{m-1}x_i;q)_{\infty}}
{(zt^{m-1}x_i;q)_{\infty}}
\prod_{i=1}^{n-m} \frac{(azt^{n-i};q)_{\infty}}{(zt^{n-i};q)_{\infty}}
\]
for the $\gsl_3$ series of Type I, and
\[
{_1\Phi_0}\hyp{a}{\text{--}}{q,t;x,y}
=\prod_{i=1}^m \frac{(azq^{-1}t^{m-1}x_i;q)_{\infty}}
{(zq^{-1}t^{m-1}x_i;q)_{\infty}}
\prod_{i=1}^{n-m} \frac{(azt^{n-i};q)_{\infty}}{(zt^{n-i};q)_{\infty}}
\]
for the $\gsl_3$ series of Type II.
\end{theorem}

If we assume $m=n$ then we may first invoke the symmetry of
Corollary~\ref{C2} to find a second pair of $q$-binomial theorems.
\begin{theorem}[(Second $\gsl_3$ $q$-binomial theorem)]
For $x=z(1,t,\dots,t^{n-1})$ and $y=(y_1,\dots,y_n)$ we have
\[
{_1\Phi_0}\hyp{a}{\text{--}}{q,t;x,y}
=\prod_{i=1}^n \frac{(azt^{n-1}y_i;q)_{\infty}}
{(zt^{n-1}y_i;q)_{\infty}}
\]
for the $\gsl_3$ series of type I, and
\[
{_1\Phi_0}\hyp{a}{\text{--}}{q,t;x,y}
=\prod_{i=1}^m \frac{(azq^{-1}t^{n-1}y_i;q)_{\infty}}
{(zq^{-1}t^{n-1}y_i;q)_{\infty}}
\]
for the $\gsl_3$ series of type II.
\end{theorem}

Using further results for $\gsl_2$ Kaneko--Macdonald series 
many more identities for $\gsl_3$ series may be proved,
such as $q$-Gauss sums, $q$-Saalsch\"utz sums, etc.
Below we restrict ourselves to just one further applications
in the form of an $\gsl_3$ analogue of Heine's
$q$-Euler transformation.
\begin{proposition}\label{propEuler}
Let $\sigma$ be fixed as in \eqref{sigma}, and
let $x=(x_1,\dots,x_m)$ and $y=z(1,t,\dots,t^{n-1})$.
Then
\[
{_2\Phi_1}\hyp{a,b}{c}{q,t;x,y}
={_2\Phi_1}\hyp{c/a,c/b}{c}{q,t;x,aby/c} 
\prod_{i=1}^m\frac{(abzq^{-\sigma}t^{m-1} x_i/c;q)_{\infty}}
{(zq^{-\sigma}t^{m-1} x_i;q)_{\infty}}
\prod_{i=1}^{n-m}\frac{(abzt^{n-i}/c;q)_{\infty}}{(zt^{n-i};q)_{\infty}}.
\]
\end{proposition}
For $b=c$ the $_2\Phi_1$ on the right is $1$ by Lemma~\ref{lemma1}
and we recover the $q$-binomial theorem of Theorem~\ref{FV}.

\begin{proof}[Proof of Proposition~\ref{propEuler}]
According to \eqref{C1}
\begin{equation*}
u_{0;z}^{(n)}\biggl(
{_2\Phi_1}\hyp{a,b}{c}{q,t;x,y}\biggr)
={_2\Phi_1}\hyp{a,b}{c}{q,t;zt^mX}.
\end{equation*}
In \cite[Proposition 3.1]{BF99} Baker and Forrester proved that
\begin{equation*}
{_2\Phi_1}\hyp{a,b}{c}{q,t;x}
={_2\Phi_1}\hyp{c/a,c/b}{c}{q,t;\frac{abx}{c}}
\prod_{i=1}^n\frac{(abx_i/c;q)_{\infty}}{(x_i;q)_{\infty}}
\end{equation*}
so that we get
\begin{equation*}
u_{0;z}^{(n)}\biggl(
{_2\Phi_1}\hyp{a,b}{c}{q,t;x,y}\biggr)
={_2\Phi_1}\hyp{c/a,c/b}{c}{q,t;\frac{abzt^m X}{c}}
\prod_{i=1}^n\frac{(abzt^m X_i/c;q)_{\infty}}{(zt^m X_i;q)_{\infty}}.
\end{equation*}
Again using \eqref{C1eq} gives
\[
u_{0;z}^{(n)}\biggl(
{_2\Phi_1}\hyp{a,b}{c}{q,t;x,y}\biggr)
=u_{0;abz/c}^{(n)}\biggl( {_2\Phi_1}\hyp{c/a,c/b}{c}{q,t;x,y} \biggr)  
\prod_{i=1}^n\frac{(abzt^m X_i/c;q)_{\infty}}{(zt^m X_i;q)_{\infty}}.
\]
Eliminating $X_i$ completes the proof.
\end{proof}

\medskip

Theorem~\ref{FV} may be viewed as a $q,t,x$-analogue 
of a result of Tarasov and Varchenko, stated in \cite[Theorem 2.3]{TV03} 
as a $\gsl_3$ discrete exponential Selberg integral.
To obtain the Tarasov--Varchenko result we take $t=q^{\gamma}$
and $a=q^{\beta+\gamma(n-1)}$ in the theorem, and let $q$ tend to $1^{-}$.
A standard computation using \eqref{ccp} and \eqref{PS} then leads to
\begin{align}\label{lim}
\sum_{\la,\mu} & z^{\abs{\mu}} v_{\la\mu}(\gamma)\,
\frac{P_{\la}^{(1/\gamma)}(x)}{P_{\la}^{(1/\gamma)}(1^m)}
\prod_{i=1}^n
\frac{\Gamma(\beta+\tilde{\mu}_i)}{\Gamma(1+\tilde{\mu}_i)}
\prod_{i=1}^m\prod_{j=1}^n
\frac{\Gamma(1-\gamma+\tilde{\la}_i-\tilde{\mu}_j)}
{\Gamma(1+\tilde{\la}_i-\tilde{\mu}_j)}  \\
&\times\prod_{1\leq i<j\leq n}
\frac{(\tilde{\mu}_i-\tilde{\mu}_j) 
\Gamma(\gamma+\tilde{\mu}_i-\tilde{\mu}_j)}
{\Gamma(1-\gamma+\tilde{\mu}_i-\tilde{\mu}_j)}
\prod_{1\leq i<j\leq m}
\frac{(\tilde{\la}_i-\tilde{\la}_j)
\Gamma(\gamma+\tilde{\la}_i-\tilde{\la}_j)}
{\Gamma(1-\gamma+\tilde{\la}_i-\tilde{\la}_j)} \notag \\
&=(1-z)^{-(\beta+(n-1)\gamma)(n-m)}
\prod_{i=1}^m(1-zx_i)^{-\beta-\gamma(n-1)}(1-x_i)^{-\gamma(m-n-1)} 
\notag \\
&\qquad\times
\prod_{i=1}^n 
\frac{\Gamma(i\gamma)\Gamma(\beta+\gamma(i-1))}{\Gamma(\gamma)}
\prod_{i=1}^m 
\frac{\Gamma(i\gamma)\Gamma(1+\gamma(i-n-1))}{\Gamma(\gamma)}.
\notag
\end{align}
Here  
\begin{equation*}
\tilde{\la}_i=\la_i+\gamma(m-i) \quad\text{and}\quad
\tilde{\mu}_i=\mu_i+\gamma(n-i),
\end{equation*}
$P_{\la}^{(\alpha)}(x)$ is the Jack polynomial:
\begin{equation*}
P_{\la}^{(\alpha)}(x)=\lim_{t\to 1} P_{\la}(x;t^{\alpha},t)
\end{equation*}
and
\begin{align*}
v_{\la\mu}(\gamma)&=\lim_{q\to 1} V_{\la\mu}(q,q^{\gamma}) \\
&=\sum_{I\subseteq [m]} (-1)^{\abs{I}} 
\prod_{\substack{i\in I \\ j\in \bar{I}}}
\frac{\tilde{\la}_j-\tilde{\la}_i+\gamma}{\tilde{\la}_j-\tilde{\la}_i}
\prod_{i\in I}\prod_{j=1}^n
\frac{\tilde{\la}_i-\tilde{\mu}_j}{\tilde{\la}_i-\tilde{\mu}_j-\gamma}.
\end{align*}
Taking $x=(w^m)$ and using the homogeneity of the Jack polynomials
(so that $P_{\mu}^{(1/\gamma)}(w^m)=
w^{\abs{\mu}}P_{\mu}^{(1/\gamma)}(1^m)$)
results in the Tarasov--Varchenko identity.
To make the correspondence exact we need to recall the
difference in normalization exhibited in \eqref{qdeformed},
and the fact that
\begin{equation*}
\prod_{i=1}^m \Gamma(1+\gamma(i-n-1))=
\frac{(-\gamma)^m n!}{(n-m)!}
\prod_{i=1}^m \Gamma(\gamma(i-n-1)).
\end{equation*}
It is interesting to note that Tarasov and Varchenko obtained the
$x=(w^m)$ instance of the series \eqref{lim}
as the coordinate function of the hypergeometric solution
of the $\gsl_3$ dynamical differential equation of \cite{TV02}
with values in the weight subspace $L_{\la}[\la-n\alpha_1-m\alpha_2]$,
$\la\in\Complex\Lambda_1$. \label{page}
Here $L_{\la}$ is an irreducible $\gsl_3$ highest
weight module of weight $\la$, and $\alpha_i$ and $\Lambda_i$ ($i=1,2$)
are the roots and fundamental weights of $\gsl_3$.
The existence of identities such as \eqref{lim} (with $x=(w^m)$)
and their associated integral evaluations was anticipated 
by Mukhin and Varchenko who formulated a very general conjecture 
regarding $\mathfrak{g}$ type Selberg integrals being expressible
in terms of products of gamma functions \cite[Conjecture~1]{MV00}.

By a standard limiting procedure the sum \eqref{lim} (with $x=(w^m)$
may be transformed into an integral, leading to the
$\gsl_3$ exponential Selberg integral of \cite[Theorem 3.1]{TV03}.
More generally, if we first transform Theorem~\ref{FV}
into a $q$-integral and then take the $q\to 1^{-}$ limit we
get a more general $\gsl_3$ Selberg integral, not contained
in \cite{TV03}.
More precisely, we take Theorem~\ref{FV} (for Type I series) and apply the
homomorphism $u_{\nu;w}^{(m)}$ acting on $x$. 
Thanks to \eqref{symm} and \eqref{deg} this yields
\begin{multline*}
\sum_{\la,\mu} w^{\abs{\la}}z^{\abs{\mu}} t^{n(\la)+n(\mu)}\, 
u_{\la}^{(n)}(P_{\nu})\, (a;q,t)_{\mu} \, \Omega_{\la\mu}(q,t) \,
\frac{u_0^{(n)}(P_{\la})}{c'_{\la}(q,t)}\,
\frac{u_0(P_{\mu})}{c'_{\mu}(q,t)}\\
=u_0^{(n)}(P_{\nu})\prod_{i=1}^m 
\frac{(awzt^{2m-i-1}q^{\nu_i};q)_{\infty}}
{(zt^{2m-i-1}wq^{\nu_i};q)_{\infty}}\:
\frac{(wt^{2m-n-i-1}q^{\nu_i};q)_{\infty}}{(wt^{m-i}q^{\nu_i};q)_{\infty}}
\prod_{i=1}^{n-m} \frac{(azt^{n-i};q)_{\infty}}{(zt^{n-i};q)_{\infty}}.
\end{multline*}
Next we replace $(a,w,z,t)\to (q^{(n-1)\gamma+\alpha},q^{\beta_1},
q^{\beta_2-m\gamma},q^{\gamma})$ and use the definition of the $q$-gamma
function to interpret this as an $(m+n)$-dimensional $q$-integral.
Taking the limit $q\to 1$ then yields a $\gsl_3$ Selberg integral
involving Jack polynomial.
The precise details of this essentially elementary calculation will be
given in a future paper in which more general Selberg-type integrals
will be considered.

To give the exact form of the integral we need to borrow some
notation from \cite{TV03}.
Let $M$ be a map
\begin{equation*}
M:\{1,\dots,m\}\to \{1,\dots,n\}
\end{equation*}
such that 
\begin{equation*}
M(i)\leq M(i+1)
\end{equation*}
and
\begin{equation*}
1\leq M(i)\leq n-m+i.
\end{equation*}
It is easily seen that there are exactly
\begin{equation*}
\frac{n-m+1}{n+1}\binom{m+n}{m}
\end{equation*}
admissible maps $M$.

Let $D^{m,n}[0,1]\subseteq [0,1]^{m+n}$ be defined as the set of points
\[
P=(x_1,\dots,x_m,y_1,\dots,y_n)
\]
such that
\begin{equation}\label{poset}
\begin{array}{ccccccccccccl}
&&&&0&\LEQ & x_1&\LEQ &x_2 &\LEQ&\dots&\LEQ&x_m \\[-1mm]
&&&&&&\text{\begin{turn}{270}$\leq$\end{turn}}&&
\text{\begin{turn}{270}$\leq$\end{turn}} &&&&
\text{\begin{turn}{270}$\leq$\end{turn}}  \\[3mm]
0&\LEQ& y_1&\LEQ&\cdots&\LEQ &y_{n-m+1}&\LEQ &y_{n-m+2}&\LEQ&\dots&\LEQ
&y_n.
\end{array}
\end{equation}
The $x$ as well as the $y$ coordinates $P\in D^{m,n}[0,1]$
are totally ordered, but only a partial order exists between the
$x_i$ and the $y_j$.
We now write $D^{m,n}[0,1]$ as a chain:
\[D^{m,n}[0,1]=\sum_M D_M^{m,n}[0,1],\]
where $D_M^{m,n}[0,1]\subseteq D^{m,n}[0,1]$ is defined
by points $P$ endowed with a total ordering among its coordinates,
by supplementing \eqref{poset} with
\begin{equation*}
y_{M(i)-1}\leq x_i\leq y_{M_s(i)}
\quad \text{for~~$1\leq i\leq m$},
\end{equation*}
where $y_0:=0$.
We further define the chain
\begin{equation}\label{chain2}
C^{m,n}_{\gamma}[0,1]=\sum_M
F_M^{m,n}(\gamma) D_M^{m,n}[0,1],
\end{equation}
where
\begin{equation*}
F_M^{m,n}(\gamma)=
\prod_{i=1}^m \frac{\sin\bigl(\pi(i+n-m-M(i)+1)\gamma\bigr)}
{\sin\bigl(\pi(i+n-m)\gamma\bigr)}.
\end{equation*}
Up to a trivial transformation (corresponding to the
variable change (5.15)) the above chains coincide 
with those of \cite{TV03}.

Finally introducing the Pochhammer symbol
\[
(a)_N=a(a+1)\cdots(a+N-1)
\]
and recalling the definition \eqref{gdef}
we are in a position to state the integral analogue of Theorem~\ref{FV}.

\begin{corollary}[($\gsl_3$ Selberg integral)]\label{CSI}
Let $\nu$ be a partition of at most $m$ parts. Then
\begin{align*}
&\Int_{C^{m,n}_{\gamma}[0,1]}P_{\nu}^{(1/\gamma)}(x) \,
h(x,y)\,
\prod_{i=1}^mx_i^{\beta_1-1} 
\prod_{i=1}^n (1-y_i)^{\alpha-1}y_i^{\beta_2-1} 
\abs{\Delta(x)}^{2\gamma}\; \abs{\Delta(y)}^{2\gamma}
\prod_{i=1}^m\prod_{j=1}^n \abs{x_i-y_j}^{-\gamma}
\; \dup x\,\dup y \\
&\qquad =\prod_{1\leq i<j\leq m}\frac{((j-i+1)\gamma)_{\nu_i-\nu_j}}
{((j-i)\gamma)_{\nu_i-\nu_j}}
\prod_{i=1}^n \frac{\Gamma(\alpha+(i-1)\gamma)\Gamma(i\gamma)}
{\Gamma(\gamma)}
\prod_{i=1}^{n-m}\frac{\Gamma(\beta_2+(i-1)\gamma)}
{\Gamma(\alpha+\beta_2+(i+n-2)\gamma)} \\
&\qquad \quad \times
\prod_{i=1}^m \frac{\Gamma(\beta_1+(m-i)\gamma+\nu_i)
\Gamma(\beta_1+\beta_2+(m-i-1)\gamma+\nu_i)
\Gamma((i-n-1)\gamma)\Gamma(i\gamma)}
{\Gamma(\beta_1+(2m-n-i-1)\gamma+\nu_i)
\Gamma(\alpha+\beta_1+\beta_2+(m+n-i-2)\gamma+\nu_i)
\Gamma(\gamma)},
\end{align*}
where
\begin{gather*}
\Re(\alpha)>0,~\Re(\beta_1)>0,~\Re(\beta_2)>0 \\[2mm]
-\min\Bigl\{
\frac{1}{n},\frac{\Re(\alpha)}{n-1},\frac{\Re(\beta_1)}{m-1},
\frac{\Re(\beta_2)}{n-m-1},\frac{\Re(\beta_1+\beta_2)}{m-2}\Bigr\}
<\Re(\gamma)<0.
\end{gather*}
\end{corollary}
The conditions on $\alpha,\beta_1,\beta_2$ and $\gamma$
(which are only sharp when $\nu=0$) are
valid for generic $n$ and $m$ and need small modifications when
$m=0,1$ or $m=n$. 
The conditions are correct for $n=1$, $m=2$
or $n=m+1$ provided $1/0$ is interpreted as $+\infty$.
Conditions that are sharp follow by demanding that
the arguments of gamma functions appearing in the numerator on the right
have positive real part. We also note that without loss of generality
one may assume that $\nu$ has at most $m-1$ parts, since
\[
P_{(\nu_1,\dots,\nu_m)}(x)=(x_1\cdots x_m)^{\nu_n}
P_{(\nu_1-\nu_m,\dots,\nu_{m-1}-\nu_m,0)}(x)
\]
so that $\nu_n$ may be eliminated by a rescaling of $\beta_1$.

For $m=0$ Corollary~\ref{CSI} is the Selberg integral \eqref{SelbergInt}
up to some trivial changes. Indeed for $m=0$ we get, after replacing 
$\beta_2$ by $\beta$,
\[
\Int_{0\leq y_1\leq \dots\leq y_n\leq 1}
\abs{\Delta(y)}^{2\gamma}
\prod_{i=1}^n (1-y_i)^{\alpha-1}y_i^{\beta-1}\: \dup y  
=\prod_{i=1}^n \frac{\Gamma(\alpha+(i-1)\gamma)
\Gamma(\beta+(i-1)\gamma)\Gamma(i\gamma)}
{\Gamma(\alpha+\beta+(i+n-2)\gamma)\Gamma(\gamma)}.
\]
Since the integrand is symmetric in $y$ and
\[
\prod_{i=1}^n \frac{\Gamma(i\gamma)}{\Gamma(\gamma)}=
\frac{1}{n!}
\prod_{i=1}^n \frac{\Gamma(i\gamma+1)}{\Gamma(\gamma+1)},
\]
this yields \eqref{SelbergInt} with $\alpha$ and $\beta$ interchanged.
(Alternatively one may replace $y_i\to 1-y_i$ for all $1\leq i\leq m$
instead of replacing $\alpha\leftrightarrow \beta$.)

When $\nu=0$ all reference to the Jack polynomial $P_{\nu}^{1/\gamma}(x)$
disappears from Corollary~\ref{CSI} and we obtain the Tarasov--Varchenko 
integral \eqref{integral}.
To make the connection with the integral of \cite{TV03} precise 
one needs to replace 
\begin{align}\label{xyst}
x_i&\to 1-s_i,& y_i&\to 1-t_i, \\ 
n&\to k_1,& m&\to k_2, \notag \\
\alpha&\to\alpha+1,& \beta_1&\leftrightarrow\beta_2 \notag 
\end{align}
and observe that
\[
h(1-s,1-t)=(-1)^{k_2}\,\tilde{h}_{k_1,k_2,k_2}(t;s)
\prod_{i=1}^{k_1} t_i^{-1},
\]
where $\tilde{h}_{l_1,l_2,m}(t;s)$ is the function defined in Section~5
of \cite{TV03}.
Then correcting a factor $(-1)^{k_2}$ missing in \cite{TV03}
one obtains the integral
\[
\tilde{J}_{k_1,k_2,k_2}(\alpha,\beta_1,\beta_2,\gamma)
\]
given by the final two equations of that paper.

\medskip

For $\nu=(1^r)$ the Jack polynomial simplifies to
the elementary symmetric function:
\[
P_{(1^r)}(x)=e_r(x)=\sum_{1\leq i_1<i_2<\dots<i_r\leq m}
x_{i_1}\cdots x_{i_r}
\]
and Corollary~\ref{CSI} yields an $\gsl_3$ version of Aomoto's
integral \cite{Aomoto95}.
\begin{corollary}
For $0\leq r\leq m$ 
\begin{align*}
&\Int_{C^{m,n}_{\gamma}[0,1]}e_r(x) \, h(x,y)\,
\prod_{i=1}^mx_i^{\beta_1-1} 
\prod_{i=1}^n (1-y_i)^{\alpha-1}y_i^{\beta_2-1} 
\abs{\Delta(x)}^{2\gamma}\; \abs{\Delta(y)}^{2\gamma}
\prod_{i=1}^m\prod_{j=1}^n \abs{x_i-y_j}^{-\gamma}
\; \dup x\,\dup y \\
&\qquad =
\binom{m}{r}
\prod_{i=1}^n \frac{\Gamma(\alpha+(i-1)\gamma)\Gamma(i\gamma)}
{\Gamma(\gamma)}
\prod_{i=1}^{n-m}\frac{\Gamma(\beta_2+(i-1)\gamma)}
{\Gamma(\alpha+\beta_2+(i+n-2)\gamma)} \\
&\qquad \quad \times
\prod_{i=1}^m\biggl(\frac{\Gamma(\beta_1+(m-i)\gamma+\chi(i\leq r))}
{\Gamma(\beta_1+(2m-n-i-1)\gamma+\chi(i\leq r))} \\
&\qquad \qquad \qquad \times
\frac{
\Gamma(\beta_1+\beta_2+(m-i-1)\gamma+\chi(i\leq r))
\Gamma((i-n-1)\gamma)\Gamma(i\gamma)}
{\Gamma(\alpha+\beta_1+\beta_2+(m+n-i-2)\gamma+\chi(i\leq r))
\Gamma(\gamma)}\biggr), 
\end{align*}
where
\begin{gather*}
\Re(\alpha)>0,~\Re(\beta_1)>0,~\Re(\beta_2)>0 \\[2mm]
-\min\Bigl\{
\frac{1}{n},\frac{\Re(\alpha)}{n-1},\frac{\Re(\beta_1)}{m-1},
\frac{\Re(\beta_2)}{n-m-1},\frac{\Re(\beta_1+\beta_2)}{m-2}\Bigr\}
<\Re(\gamma)<0
\end{gather*}
and $\chi(\textup{true})=1$, $\chi(\textup{false})=0$.
\end{corollary}
The comments made immediately after Corollary~\ref{CSI} still apply.

\subsection*{Acknowledgements}
I am much indebted to Michael Schlosser for very helpful
discussions and for pointing out the significance of
\cite{KN99} for our work.

\bibliographystyle{amsplain}

\begin{thebibliography}{99}

\bibitem{AAR99}
G. E. Andrews, R. Askey and R. Roy,
\textit{Special functions},
Encyclopedia of Mathematics and its Applications, Vol.~71,
(Cambridge University Press, Cambridge, 1999).

\bibitem{Aomoto95}
K. Aomoto,
\textit{Connection formulas of the $q$-analog de Rham cohomology},
in \textit{Functional Analysis on the Eve of the 21st Century}, Vol.~1,
pp. 1--12, S. Gindikin \textit{et al.} eds.,
Prog. in Math. 131 (Birkhauser, Boston, MA, 1995).

\bibitem{BF99}
T. H. Baker and P. J. Forrester,
\textit{Transformation formulas for multivariable basic 
hypergeometric series},
Methods Appl. Anal. \textbf{6} (1999), 147--164.

\bibitem{BS98}
G. Bhatnagar and M. Schlosser,
\textit{$C_n$ and $D_n$ very-well-poised $_{10}\phi_9$ transformations},
Constr. Approx. \textbf{14} (1998), 531--567.

\bibitem{Bressoud99}
D. M. Bressoud,
\textit{Proofs and confirmations --- The story of the alternating 
sign matrix conjecture},
(Cambridge University Press, Cambridge, 1999).

\bibitem{GR04}
G. Gasper and M. Rahman,
\textit{Basic Hypergeometric Series}, 
Encyclopedia of Mathematics and its Applications, Vol.~35,
second edition,
(Cambridge University Press, Cambridge, 2004).

\bibitem{GK96}
R. A. Gustafson and C. Krattenthaler,
\textit{Heine transformations for a new kind of basic hypergeometric series in $U(n)$}, 
J. Comput. Math. Appl. \textbf{68} (1996), 151--158.

\bibitem{Izergin87}
A. G. Izergin,
\textit{Partition function of a six-vertex model in a finite volume},
Dokl. Akad. Nauk SSSR \textbf{297} (1987), 331--333.

\bibitem{Kaneko96}
J. Kaneko,
\textit{$q$-Selberg integrals and Macdonald polynomials},
Ann. Sci. \'Ecole Norm. Sup. \textbf{29} (1996), 583--637.

\bibitem{KN99}
A. N. Kirillov and M. Noumi,
\textit{$q$-Difference raising operators for Macdonald polynomials
and the integrality of transition coefficients} in
\textit{Algebraic methods and $q$-special functions}, pp. 227--243,
CRM Proc. Lecture Notes 22,
(Amer. Math. Soc., Providence, RI, 1999). 

\bibitem{Korepin82}
V. E. Korepin,
\textit{Calculation of norms of Bethe wave functions},
Comm. Math. Phys. \textbf{86} (1982), 391--418.

\bibitem{Krattenthaler99}
C. Krattenthaler,
\textit{Advanced determinant calculus},
S\'em. Lothar. Combin. \textbf{42} (1999), Art. B42q, 41pp.

\bibitem{Kuperberg96}
G. Kuperberg,
\textit{Another proof of the alternating-sign matrix conjecture},
Internat. Math. Res. Notices (1996), 139--150.

\bibitem{Lascoux99}
A. Lascoux,
\textit{Square-ice enumeration},
S\'em. Lothar. Combin. \textbf{42} (1999), Art. B42p, 15pp.

\bibitem{Macdonald82}
I. G. Macdonald,
\textit{Some conjectures for root systems},
SIAM J. Math. Anal. \textbf{13} (1982), 988--1007.

\bibitem{Macdonald95}
I. G. Macdonald,
\textit{Symmetric functions and Hall polynomials},
second edition,
(Oxford University Press, New York, 1995).

\bibitem{Macdonald}
I. G. Macdonald,
\textit{Hypergeometric series II}, unpublished manuscript.

\bibitem{Milne85}
S. C. Milne,
\textit{
An elementary proof of the Macdonald identities for $A^{(1)}_l$},
Adv. in Math. \textbf{57} (1985), 34--70.

\bibitem{Milne92}
S. C. Milne,
\textit{Summation theorems for basic hypergeometric series of
Schur function argument},
in \textit{Progress in approximation theory}, pp. 51--77,
A.~A.~Gonchar and E.~B.~Saff eds.,
Springer Ser. Comput. Math. 19 (Springer, New York, 1992).

\bibitem{Milne97}
S. C. Milne,
\textit{Balanced $_3\phi_2$ summation theorems for $U(n)$
basic hypergeometric series},
Adv. Math. \textbf{131} (1997), 93--187.

\bibitem{ML95}
S. C. Milne and G. M. Lilly,
\textit{Consequences of the $A_{\ell}$ and $C_{\ell}$ Bailey transform and
Bailey lemma},
Discrete Math. \textbf{139} (1995), 319--346.

\bibitem{MV00}
E. Mukhin and A. Varchenko,
\textit{Remarks on critical points of phase functions and norms of 
Bethe vectors},
Adv. Stud. Pure Math. \textbf{27} (2000), 239--246.

\bibitem{Selberg44}
A. Selberg,
\textit{Bemerkninger om et multipelt integral},
Norske Mat. Tidsskr. \textbf{26} (1944), 71--78.

\bibitem{TV02}
V. Tarasov and A. Varchenko, 
\textit{Duality for Knizhnik--Zamolodchikov and dynamical equations},
Acta Appl. Math. \textbf{73} (2002), 141--154.

\bibitem{TV03}
V. Tarasov and A. Varchenko, 
\textit{Selberg-type integrals associated with $\gsl_3$},
Lett. Math. Phys. \textbf{65} (2003), 173--185.

\bibitem{Warnaar}
S. O. Warnaar,
\textit{$q$-Selberg integrals and Macdonald polynomials},
Ramanujan J. \textbf{10} (2005), 237--268.

\bibitem{Warnaarb}
S. O. Warnaar,
\textit{A Selberg integral for the Lie algebra A$_n$},
in preparation.

\end{thebibliography}

\end{document}